\documentclass[12pt]{article}
\usepackage{amssymb}
\usepackage{amsmath}
\usepackage{amsfonts}
\usepackage[T1]{fontenc}
\usepackage{graphicx}
\usepackage{epic}

\newcommand{\R}{\mathbb{R}}

\newcommand{\dee}{\mathop{\! \, \mathrm{d} \!}\nolimits}

\allowdisplaybreaks

\begin{document}

\title{Bohr-Sommerfeld-Heisenberg Theory \\
in Geometric Quantization}
\author{Richard Cushman and J\k{e}drzej \'{S}niatycki\thanks{%
Department of Mathematics and Statistics, University of Calgary, Calgary,
Alberta, T2N 1N4 Canada}}
\date{}
\maketitle

\begin{abstract}
In the framework of geometric quantization we extend the Bohr-Sommerfeld
rules to a full quantum theory which resembles Heisenberg's matrix theory.
This extension is possible because Bohr-Sommer-feld rules not only provide
an orthogonal basis in the space of quantum states, but also give a lattice
structure to this basis. This permits the definition of appropriate shifting
operators. As examples, we discuss the $1$-dimensional harmonic oscillator
and the coadjoint orbits of the rotation group.
\end{abstract}


\section{Introduction}

\label{sec1}

\bigskip

The desire to understand the energy spectrum of completely integrable
Hamiltonian systems lead to Bohr-Sommerfeld theory, also called old quantum
theory. Bohr \cite{bohr} explained Planck's hypothesis by the spectrum of
the harmonic oscillator, which he obtained using his quantum conditions.
Sommerfeld \cite{sommerfeld} extended Bohr's quantization rules to a system
with Hamiltonian 
\begin{equation}
H=\sqrt{\boldsymbol{p}^{2}+m^{2}}+%
\mbox{$\frac{{\scriptstyle
k}}{{\scriptstyle r}}$}.  \label{I1}
\end{equation}%
The first term in (\ref{I1}) is the relativistic expression for the kinetic
energy of a particle with mass $m$ and momentum $\boldsymbol{p}$ and the
second term is the potential energy of a charged particle in the electric
field produced by a stationary charged particle at the origin. With an
appropriate choice of the parameters $m$ and $k$, we can use the Hamiltonian 
$H$ as an approximation to the energy of an electron in the hydrogen atom in
the limit of infinite mass of its nucleus. The energy spectrum of the
hydrogen atom obtained by Sommerfeld agreed exactly with the observed
spectrum.\footnote{%
It is remarkable that the energy spectrum obtained by Sommerfeld agrees
exactly with the energy spectrum obtained by solving the Dirac equation for
an electron in the same electric field \cite{davydov}. Even more puzzling is
the fact that a modification of Bohr-Sommerfeld conditions by a term $\frac{1%
}{2}\hbar $ gives rise to the energy spectrum for a $\pi $ meson in the same
electric field, whixh can be obtained by solving the Klein-Gordon equation. 
\cite{sniatycki75}} The Bohr-Sommerfeld theory was applied with varying
success to other systems. \medskip

The problem with Bohr-Sommerfeld theory is that it gives only the joint
spectrum of energy and angular momentum. It does not provide a way to
discuss the probability of transition between states. The next stage in the
development of understanding of the nature of quantum physics was provided
by the matrix theory of Heisenberg \cite{heisenberg} and the wave theory of
Schr\"{o}dinger \cite{schrodinger}. Heisenberg postulated that dynamical
variables were not functions on the phase space of the system but matrices
in some vector space, possibly infinite dimensional. One can infer that
Heisenberg's matrices are linear transformations in the space of physical
states relative to a basis provided by Bohr-Sommerfeld joint eigenstates of
energy and angular momentum. Heisenberg's approach was further developed by
Born and Jordan \cite{born-jordan}, who used it to study various physical
systems. For Schr\"{o}dinger, physical states were described by complex
valued functions on the configuration space of the system; while dynamical
variables were represented by differential operators obtained from classical
dynamical variables by replacing the momentum by $i/\hbar $ times the
operator of differentiation with respect to the position variable. Here $%
\hbar $ is Planck's constant divided by $2\pi $. Dirac \cite{dirac} showed
that the theories of Heisenberg and Schr\"{o}dinger are equivalent. Since
then, the Schr\"{o}dinger equation has become the computational basis of
quantum mechanics. Heisenberg's theory is discussed in works mainly of
historical interest \cite{mehra-rechenberg}. \medskip

At present, quantization of a completely integrable Hamiltonian system is
discussed in the Schr\"{o}dinger framework. The energy spectra obtained
there tend to their Bohr-Sommerfeld counterparts as $\hbar \rightarrow 0$ 
\cite{vungoc}. \medskip

The aim of this paper is to find a place for the Heisenberg matrix formalism
within the framework of geometric quantization. A completely integrable
Hamiltonian system defines a singular real polarization $F=D\otimes \mathbb{C%
}$ of the phase space of the underlying classical system, which is a
symplectic manifold $(P,\omega )$. We denote by $\mathcal{S}_{F}^{\infty
}(L) $ the space of sections of the prequantization line bundle $L$ over $P$
that are covariantly constant along $F$. The Bohr-Sommerfeld conditions
identify those leaves of $F$ that admit lifts to covariantly constant
sections of $L$ \cite{sniatycki}. The space $\mathcal{S}_{F}^{\infty }(L)$
consists of generalized sections (distribution sections) of $L$ that are
supported on unions of Bohr-Sommerfeld leaves. For each Bohr-Sommerfeld
leaf, we can choose a section supported on that leaf. In this way, we obtain
a basis of complex vector space $\mathcal{S}_{F}^{\infty }(L)$. We may
choose a scalar product on $\mathcal{S}_{F}^{\infty }(L)$ so that this basis
is orthogonal. Let $\mathfrak{H}$ denotes the space of sections of $\mathcal{%
S}_{F}^{\infty }(L)$ that are normalizable with respect to this scalar
product. \medskip

It is natural to choose a scalar product in $\mathcal{S}_{F}^{\infty }(L)$
so that sections supported on different Bohr-Sommerfeld leaves are mutually
orthogonal. However, the classical theory does not suggest how to normalize
any section supported on a single Bohr-Sommerfeld leaf, as this
normalization is quite arbitrary. Nevertheless, the projection of our
orthonormal basis to the complex projective space of one dimensional
subspaces of $\mathfrak{H}$ is well defined by the classical data. \medskip

Our next step is the observation that the basis of $\mathfrak{H}$ given by
the Bohr-Sommerfeld conditions has a natural structure of a local lattice.
This observation lead Cushman and Duistermaat \cite{cushman-duistermaat} to
the notion of a quantum monodromy. First suppose that our basis of $%
\mathfrak{H}$ is a global lattice. Then, associated to each generator of the
lattice, there is a well defined shifting operator defined by assigning to
each vector of our basis the adjacent vector in the lattice. Next we compute
commutation relations among the shifting operators and the operators that
are diagonal in our basis. Furthermore we look for functions on $P$, which
satisfy the Poisson bracket relations corresponding to the commutation
relations of shifting operators. We define quantization of these functions
to be the corresponding shifting operators. As usual in quantization, we are
making a choice, but we know exactly the arbitrariness involved in this
choice. We apply the above this procedure to harmonic oscillator with one
degree of freedom and to coadjoint orbits of the rotation group $%
\mathop{\mathrm{SO}}\nolimits (3)$. In \cite{cushman-sniatycki12a} we treat
the harmonic oscillator in two degrees of freedom. In these examples, our
theory gives the usual results. \medskip

Now suppose that the Bohr-Sommerfeld basis is only a local lattice. Let $U$
be the maximal open dense subset of $P$ on which the singular real
polarization $F$ of $(P,\omega )$ is regular. The space $V=U/D$ of the
regular tori in $P$ is a quotient manifold of $U$ and the projection map $%
\pi :U\rightarrow V$ is a locally trivial $2$-torus bundle. In the universal
covering space $\widetilde{V}$ of $V,$ we repeat the construction for the
case of a global lattice. Quantization of functions on $\widetilde{V}$ that
are pull-backs of functions from $V$ gives quantum operators in our Hilbert
space $\mathfrak{H}$. We illustrate our approach in another paper \cite%
{cushman-sniatycki12c} by quantizing the spherical pendulum.

\section{Elements of geometric quantization}


In this section, we review the elements of geometric quantization. Our
notation differs from other authors in that our symplectic form is the \emph{%
negative} of the symplectic form used in \cite{kostant}. This is the reason
for the appearance of a negative sign in various formulae. In particular, if 
$(P,\omega )$ is a symplectic manifold, then the Hamiltonian vector field $%
X_{f}$ of corresponding to the Hamiltonian function $f\in C^{\infty }(P)$
satisfies $X_{f} \mbox{$\, \rule{8pt}{.5pt}\rule{.5pt}{6pt}\, \, $} \omega =
-\mathop{\! \, \mathrm{d} \!}\nolimits f$ and the Poisson bracket of $f_{1},
f_{2}\in C^{\infty }(P)$ is given by $\{ f_{1},f_{2}\}= X_{f_{2}}f_{1}$.
These conventions are commonly used in theoretical physics, see \cite%
{sternberg}, \cite{sniatycki}, and \cite{woodhouse}.

\subsection{Prequantization}


We now discuss prequantization. \medskip

Let $\lambda :L\rightarrow P$ be a complex line bundle with a connection and
a connection invariant Hermitian inner product $\langle \cdot \mid \cdot
\rangle $. A connection on $L$ is given by the covariant derivative operator 
$\nabla $, which associates to each section $\sigma $ of $L$ and each smooth
vector field $X$ on $P$ a section $\nabla_{X}\sigma $ of $L$ so that for
each $f\in C^{\infty}(P) $, 
\begin{equation*}
\nabla _{X}(f\sigma )=X(f)\sigma +f\nabla _{X}\sigma \mspace{10mu} \text{and}%
\mspace{10mu}\nabla_{fX}\sigma =f\nabla _{X}\sigma \text{.}
\end{equation*}
For every section $\sigma $ of $L$, $f\in C^{\infty }(P)$ and every smooth
vector field $X,$ $X^{\prime }$ on $P$, the curvature $K(X,X^{\prime }) =
[\nabla _{X}, \nabla _{X^{\prime }}] - {\nabla }_{[X,X^{\prime }]}$ of $%
\nabla $ satisfies $K(X^{\prime }, X) = -K(X, X^{\prime })$ and $%
K(X,X^{\prime })(f\, \sigma ) = f \, K(X,X^{\prime })\sigma $. Hence there
is a $2$-form $\alpha $ on $P$ such that 
\begin{equation}
K(X,X^{\prime })\sigma = 2\pi i\, \alpha (X,X^{\prime })\sigma \text{.}
\label{7curvature}
\end{equation}%
The form $\alpha $ is the pull-back by the section $\sigma $ of the
curvature form of the connection $\nabla $. The Hermitian form $\langle
\cdot \mid \cdot \rangle $ on $L$ is connection invariant if, for\ every
pair of sections $\sigma _{1},\sigma _{2}$ of $L$ and every smooth vector
field $X$ on $P$ we have 
\begin{equation*}
X(\langle \sigma _{1}\mid \sigma _{2}\rangle )=\langle \nabla _{X}\sigma
_{1}\mid \sigma _{2}\rangle +\langle \sigma _{1}\mid \nabla _{X}\sigma
_{2}\rangle \text{.}
\end{equation*}

Quantization of a mechanical system is defined in terms of an additional
free parameter $\hbar $. In quantum mechanics, $\hbar $ is the value of
Planck's constant divided by $2\pi $. However, in the quasi-classical
approximation one considers limits of various expressions as $\hbar
\rightarrow 0$. \medskip

The line bundle $L$ over $P$ with a connection $\nabla $ and a connection
invariant Hermitian form on $L$ is a prequantization line bundle of $%
(P,\omega )$ if the following prequantization condition is satisfied 
\begin{equation}
K(X,X^{\prime })\sigma =%
\mbox{$\frac{{\scriptstyle i}}{{\scriptstyle
\hbar}}$}\,\omega (X,X^{\prime })\sigma   \label{PCondition1}
\end{equation}%
for every smooth vector field $X,X^{\prime }$ on $P$ and each section $%
\sigma $ of $L$. The prequantization condition (\ref{PCondition1}) requires
that the de Rham cohomology class $[(2\pi \hbar )^{-1}\omega ]$ on $P$ is in 
${\mathrm{H}}^{2}(P,\mathbb{Z})$. \medskip 

Prequantization assigns to each $f\in C^{\infty }(P)$ an operator ${\mathbf{P%
}}_f$ on the space $S^{\infty }(L)$ of smooth sections of $L$ given by 
\begin{equation}
\boldsymbol{P}_{f}\sigma = i\hbar (\nabla _{X_{f}}+f)\sigma \text{.}
\label{7prequantization}
\end{equation}%
For each $f_{1},f_{2}\in C^{\infty }(P)$ and $\sigma \in S^{\infty }(L),$ we
have%
\begin{equation}
\lbrack \boldsymbol{P}_{f_{1}},\boldsymbol{P}_{f_{2}}]=-i \hbar \, 
\boldsymbol{P}_{\{f_{1},f_{2}\}}\text{.}  \label{7Poisson}
\end{equation}%
This implies that the map 
\begin{equation*}
C^{\infty }(P)\times S^{\infty }(L)\rightarrow S^{\infty
}(L):(f,\sigma)\mapsto 
\mbox{$\frac{{\scriptstyle i }}{{\scriptstyle \hbar
}}$}P_{f}\sigma
\end{equation*}
is a representation of the Lie algebra of $(C^{\infty }(P), \{ \, \, , \, \,
\})$ on $S^{\infty }(L).$ \medskip

The space $S_{0}^{\infty }(L)$ of compactly supported smooth sections of $L$
has a Hermitian inner product 
\begin{equation}
(\sigma _{1}\mid \sigma _{2})=\int_{P}\langle \sigma _{1}\mid \sigma
_{2}\rangle \omega ^{n},  \label{7ScalarProduct}
\end{equation}%
where $n=\mbox{$\frac{\scriptstyle 1}{\scriptstyle 2}\,$} \dim P$. For each $%
f\in C^{\infty }(P)$, the prequantization operator $\boldsymbol{P}_{f}$ is
symmetric with respect to the inner product (\ref{7ScalarProduct}). If the
Hamiltonian vector field $X_{f}$ of $f$ is complete, then $\boldsymbol{P}%
_{f} $ is self adjoint on the Hilbert space obtained by completing $%
S_{0}^{\infty }(L)$ with respect to the norm given by (\ref{7ScalarProduct}%
). Equation (\ref{7Poisson}) gives the usual commutation relations imposed
in quantum mechanics.

\subsection{\textbf{Polarization}}


Prequantization does not correspond to quantum theory because the
probability density $(\sigma \! \mid \! \sigma )(p)$ of localizing the state 
$\sigma $ at a point $p\in P$ fails to satisfy the Heisenberg uncertainty
principle. To avoid this difficulty we introduce the notion of a
polarization. \medskip

A complex distribution $F\subset T^{\mathbb{C}}P=\mathbb{C} \otimes TP$ on a
symplectic manifold $(P,\omega )$ is Lagrangian if, for each $p\in P$, the
restriction of the symplectic form $\omega $ to the subspace $F_{p}\subset
T_{p}^{\mathbb{C}}P$ vanishes identically and $\mathrm{rank}_{\mathbb{C}}\,F=%
\mbox{$\frac{{\scriptstyle 1}}{{\scriptstyle 2}}$}\dim P$. We denote the
complex conjugate of the distribution $F$ by $\overline{F}$. Let 
\begin{equation}
D=F\cap \overline{F}\cap TP\mspace{10mu}\text{and}\mspace{10mu} E=(F+%
\overline{F})\cap TP\text{.}  \label{eq-5distrib}
\end{equation}%
A polarization of $(P,\omega )$ is an involutive complex Lagrangian
distribution $F$ such that $D$ and $E$ are involutive distributions on $P$.
\medskip

Let $C^{\infty }(P)_{F}^{0}$ be the space of smooth complex valued functions
on $P$ that are constant along $F$, that is, 
\begin{equation*}
C^{\infty }(P)_{F}^{0}=\{f\in C^{\infty }(P)\otimes \mathbb{C}\,\mid uf=0%
\text{ for all }u\in F\}.
\end{equation*}%
We denote by $C_{F}^{\infty }(P)$ the space of smooth functions on $P$ whose
Hamiltonian vector fields preserve $F$. In other words, $f\in C_{F}^{\infty
}(P)$ if, for every $h\in C^{\infty }(P)_{F}^{0},$ the Poisson bracket $%
\{f,h\}\in C^{\infty }(P)_{F}^{0}$. If $f_{1},f_{2}\in C_{F}^{\infty }(P)$
and $h\in C^{\infty }(P)_{F}^{0}$ then the Jacobi identity implies that 
\begin{equation*}
\{\{f_{1},f_{2}\},h\}=-\{f_{2},\{f_{1},h\}\}+\{f_{1},\{f_{2},h\}\}\in
C^{\infty }(P)_{F}^{0}\text{.}
\end{equation*}%
Hence, for a strongly admissible polarization, the ring $C_{F}^{\infty }(P)$
is a Poisson subalgebra of $(C^{\infty }(P),\{\,\,,\,\,\})$. \medskip

Let $\mathcal{S}_{F}^{\infty }(L)$ denote the space of smooth sections of $L$
which are covariantly constant along $F$, namely, 
\begin{equation*}
\mathcal{S}_{F}^{\infty }(L)=\{\sigma \in \mathcal{S}^{\infty }(L)\mid
\nabla _{u}\sigma =0\mspace{5mu}\text{for all}\mspace{5mu}u\in F\}.
\end{equation*}%
For each $h\in C^{\infty }(P)_{F}^{0},$ $f\in C_{F}^{\infty }(P)$ and $%
\sigma \in \mathcal{S}_{F}^{\infty }(L)$ we have $\nabla _{X_{h}}(%
\boldsymbol{P}_{f}\sigma )=0$. Thus, for every $f\in C_{F}^{\infty }(P)$,
the prequantization operator $\boldsymbol{P}_{f}$ maps $\mathcal{S}%
_{F}^{\infty }(L)$ to itself. The quantization map $\boldsymbol{Q}$ relative
to a polarization $F$ is the restriction of the prequantization map 
\begin{equation*}
\boldsymbol{P}:C^{\infty }(P)\times \mathcal{S}^{\infty }(L)\rightarrow 
\mathcal{S}^{\infty }(L):(f,\sigma )\mapsto \boldsymbol{P}_{f}\sigma =
i\hbar (\nabla _{X_{f}}+f)\sigma
\end{equation*}%
to the domain $C_{F}^{\infty }(P)\times \mathcal{S}_{F}^{\infty }(L)\subset
C^{\infty }(P)\times \mathcal{S}^{\infty }(L)$ and the codomain $\mathcal{S}%
_{F}^{\infty }(L)\subset \mathcal{S}^{\infty }(L)$. In other words, 
\begin{equation}
\boldsymbol{Q}:C_{F}^{\infty }(P)\times \mathcal{S}_{F}^{\infty
}(L)\rightarrow \mathcal{S}_{F}^{\infty }(L):(f,\sigma )\mapsto \boldsymbol{Q%
}_{f}\sigma = i \hbar (\nabla _{X_{f}}+f)\sigma \text{.}
\label{7Quantization}
\end{equation}%
For each $f_{1},f_{2}\in C_{F}^{\infty }(P)$, the quantized operators ${%
\mathbf{Q}}_{f_{1}}$ and ${\mathbf{Q}}_{f_{2}}$ satisfy the Dirac
commutation relations 
\begin{equation}
\lbrack \boldsymbol{Q}_{f_{1}},\boldsymbol{Q}_{f_{2}}]=-i \hbar \,%
\boldsymbol{Q}_{\{f_{1},f_{2}\}}  \label{Dirac}
\end{equation}

The choice of a polarization in geometric quantization is analogous to the
choice of a complete family of commuting observables in the Dirac theory. In
this paper we choose a real polarization corresponding to a foliation of $%
(P,\omega )$ by Lagrangian tori. \medskip

In general, sections in $S_{F}^{\infty }(L)$ need not be square integrable
with respect to the inner product (\ref{7ScalarProduct}). Therefore, one may
have to introduce a new inner product in $S_{F}^{\infty }(L)$. We refer to
this step as unitarization. In the situation considered here, the choice of
the inner product will be discussed later.

\section{The Bohr-Sommerfeld conditions}


A Hamiltonian system on $(P,\omega )$ is completely integrable if it admits $%
n=\mbox{$\frac{\scriptstyle 1}{\scriptstyle 2}$} \dim P$ Poisson commuting
constants of motion $f_{1},...,f_{n}$, which are functionally independent on
an open dense subset $U$ of $P$ and the joint level sets of $f_{1},...,f_{n}$
form a singular foliation of $(P,\omega )$ by $n$-dimensional Lagrangian
tori. In other words, the span of the Hamiltonian vector fields of $%
f_{1},...,f_{n}$ defines a singular real polarization $D$ of $(P,\omega )$.
The restriction of $D$ to the open dense subset $U$ of $P$ is a regular
polarization of $(U,\omega _{\mid U})$. \medskip

Since leaves of $D$ are affine Lagrangian $n$-tori, the connection on $L$
restricted to each leaf is flat. Hence, the only obstruction to the
existence of sections of $L$ that are covariantly constant along $D$ is the
vanishing of the holonomy group. Let $T$ be an integral manifold of $D$.
From the existence of action-angle coordinates it follows that there is a
neighbourhood $W$ of $T$ in $P$ such that $\omega $ restricted to $W$ is
exact, that is $\omega _{\mid W}=\mathop{\! \, \mathrm{d}
\!}\nolimits\theta _{W}$ for a $1$-form $\theta _{W}$ on $W$, see \cite[%
appendix D]{cushman-bates}. The holonomy group of $L_{\mid T}$ vanishes if
an only if, for each generator $\Gamma _{i}$ of the fundamental group of the 
$n$-torus $T,$ 
\begin{equation}
\int_{\Gamma _{i}}\theta =m_{i}h\mspace{10mu}\text{for}\mspace{5mu}i=1,...,n,
\label{Bohr-Sommerfeld}
\end{equation}%
where $m_{i}$ is an integer and $h$ is Planck's constant. For proof of this
statement see \cite{sniatycki}. Equation (\ref{Bohr-Sommerfeld}) is known as
the Bohr-Sommerfeld quantization condition. \medskip

Let $S$ be the collection of all tori satisfying the Bohr-Sommerfeld
condition. We refer to $S$ as the Bohr-Sommerfeld set of the integrable
system $(f_1, \ldots , f_n, P, \omega )$. Since the curvature form of $L$ is
symplectic, it follows that the complement of $S $ is open in $P$. Hence,
the representation space $\mathfrak{H}$ of geometric quantization of an
integrable system consists of distribution sections of $L$ supported on the
Bohr Sommerfeld set $S$. Since these distribution sections are covariantly
constant along the distribution $D$, it follows that each $n$-torus $T\in S$
corresponds to a 1-dimensional subspace $\mathfrak{H}_{T}$ of $\mathfrak{H}$%
. We choose a inner product $(\, \, | \, \, )$ on $\mathfrak{H}$ so that the
family $\{\mathfrak{H}_{T}\, \mathop{\rule[-4pt]{.5pt}{13pt}\, }\nolimits \,
T\in S\}$ consists of mutually orthogonal subspaces. \medskip

In Bohr-Sommerfeld quantization, one assigns to each $n$-tuple of Poisson
commuting constants of motion $f=(f_{1},...,f_{n})$ on $P$ an $n$-tuple $(%
\boldsymbol{Q}_{f_{1}},\ldots ,\boldsymbol{Q}_{f_{n}})$ of commuting quantum
operators $\boldsymbol{Q}_{f_{k}}$ for $1\leq k\leq n$ such that for each $n$%
-torus $T\in S$, the corresponding $1$-dimensional space $\mathfrak{H}_{T}$
of the representation space $(\mathfrak{H},(\,\,|\,\,))$ is an eigenspace
for each $\boldsymbol{Q}_{f_{k}}$ for $1\leq k\leq n$ with eigenvalue $%
f_{k\mid T}$. For any smooth function $F\in C^{\infty }(\mathbb{R}^{n}),$
the composition $F(f_{1},...,f_{n})$ is quantizable. The operator $%
\boldsymbol{Q}_{F(f_{1},...,f_{n})}$ acts on each $\mathfrak{H}_{T}$ by
multiplication by $F(f_{1},...,f_{n})_{\mid T}$. \medskip 

The disadvantage of the Bohr-Sommerfeld theory, as described above, is that
it does not give rise to operators describing transitions between
corresponding tori in $S$. Nevertheless, it has been useful in determining
the dimension of the space of states of certain quantum systems, see \cite%
{guillemin-sternberg} and \cite{jeffrey-weitsman}. In this paper, we follow
ideas due to Heisenberg \cite{heisenberg} and Born and Jordan \cite%
{born-jordan} to obtain an extension of the Bohr-Sommerfeld theory to a full
quantum theory with a large class of quantizable functions.

\section{Shifting operators}


From the Bohr-Sommerfeld condition (\ref{Bohr-Sommerfeld}) it follows that
the Bohr-Sommer-feld set $S$ is a local lattice. In this section we assume
that there exist global action-angle variables $(A_{i},\varphi _{i})$ on $U$
such that $\omega _{\mid U}=\mathop{\! \, \mathrm{d} \!}\nolimits\,(%
\sum_{i=1}^{n}A_{i}\mathop{\! \, \mathrm{d} \!}\nolimits\varphi _{i})$. In
other words, we assume that the $1$-form $\theta =\sum_{i=1}^{n}A_{i}\,%
\mathop{\! \, \mathrm{d} \!}\nolimits\varphi _{i}$ is defined \emph{globally}
on $U.$ In this case the Bohr Sommerfeld set $S$ defines a global lattice $%
S_{U}$ on $U$. The case when the action-angle variables are only locally
defined will be studied in \cite{cushman-sniatycki12c}. \medskip

Consider a subspace ${\mathfrak{H}}_U$ of $\mathfrak{H}$ given by the direct
sum of $1$-dimensional subspaces of ${\mathfrak{H}}_U$ corresponding to $n$%
-tori \ $T\in S_U$. Since $S_U$ is a global lattice, we can label the $n$%
-tori in $S_U$ by $n$-tuples of integers $\mathbf{m}=(m_{1},...,m_{n})$,
where the index $i$ corresponds to index of the action-angle variables $%
(A_{i},\,\varphi _{i})$ and $m_{i}$ is the integer appearing in the
Bohr-Sommerfeld condition (\ref{Bohr-Sommerfeld}), which defines an $n$%
-torus in these variables. In other words, we write the Bohr-Sommerfeld
conditions in the form 
\begin{equation*}
\int_{{\Gamma }_{i}}A_{i}\, \mathop{\! \, \mathrm{d} \!}\nolimits
\varphi_{i}=m_{i}h \mspace{10mu}\text{for each}\mspace{5mu}i=1,...,n.
\end{equation*}%
Since the actions $A_{i}$ are independent of the angle variables, we can
perform the integration and obtain 
\begin{equation}
A_{i}=m_{i}\hbar \mspace{10mu}\text{for each}\mspace{5mu}i=1,...,n.
\label{Bohr-Sommerfeld1}
\end{equation}%
where $\hbar =h/2\pi $. Equation (\ref{Bohr-Sommerfeld1}) determine an $n$%
-torus $T_{\mathbf{m}}$ in $S_{U}$. Let $\boldsymbol{e}_{\mathbf{m}}$ be a
basis vector of $\mathfrak{H}_{T_{\mathbf{m}}}$ corresponding to the $n$%
-torus $T_{\mathbf{m}}$. Each $\boldsymbol{e}_{\mathbf{m}}$ is a joint
eigenvector of the commuting operators $(\boldsymbol{Q}_{A_{1}},\ldots , 
\boldsymbol{Q}_{A_{n}})$ corresponding to eigenvalues $(m_{1}\hbar ,\ldots
,m_{n}\hbar )$. The vectors $(\boldsymbol{e}_{\mathbf{m}})$ form an
orthonormal basis in ${\mathfrak{H}}_U$. Thus, 
\begin{equation}
\left( \boldsymbol{e}_{\mathbf{m}}\mid \boldsymbol{e}_{{\mathbf{m}}^{\prime
}}\right) =0\mspace{10mu}\text{if}\mspace{5mu} \mathbf{m}\neq {\mathbf{m}}%
^{\prime }.  \label{product1}
\end{equation}

For each $i=1,...,n$, introduce an operator $\boldsymbol{a}_{i}$ on ${%
\mathfrak{H}}_{U}$ such that 
\begin{equation}
\boldsymbol{a}_{i}\boldsymbol{e}%
_{(m_{1},..,,m_{i-1},m_{i},m_{i+1},...,m_{n})}=\boldsymbol{e}%
_{(m_{1},...m_{i-1},m_{i}-1,m_{i+1},...,m_{n})\text{.}}  \label{5one}
\end{equation}%
In other words, the operator $\boldsymbol{a}_{i}$ shifts the joint
eigenspace of $(\boldsymbol{Q}_{A_{1}},\ldots ,\boldsymbol{Q}_{A_{n}})$
corresponding to the eigenvalue $(m_{1}\hbar ,\ldots ,m_{n}\hbar )$ to the
joint eigenspace of $(\boldsymbol{Q}_{A_{1}},\ldots ,\boldsymbol{Q}_{A_{n}})$
corresponding to the eigenvalue $(m_{1}\hbar ,$ $\ldots ,m_{i-1}\hbar
,(m_{i}-1)\hbar ,$ $m_{i+1}\hbar ,\ldots ,m_{n}\hbar )$. Let $\boldsymbol{a}%
_{i}^{\dagger }$ be the adjoint of $\boldsymbol{a}_{i}$. Equations (\ref%
{product1}) and (\ref{5one}) yield 
\begin{equation}
\boldsymbol{a}_{i}^{\dagger }\boldsymbol{e}%
_{(m_{1},..,,m_{i-1},m_{i},m_{i+1},...,m_{n})}=\boldsymbol{e}%
_{(m_{1},...m_{i-1},m_{i}+1,m_{i+1},...,m_{n})}  \label{6}
\end{equation}%
We refer to the operators $\boldsymbol{a}_{i}$ and $\boldsymbol{a}%
_{i}^{\dagger }$ as \emph{shifting} operators.\footnote{%
In representation theory, shifting operators are called ladder operators.
The corresponding operators in quantum field theory are called the creation
and annihilation operators.} For every $i=1,...,n$ and each $\mathbf{m}$, we
have 
\begin{align}
\lbrack \boldsymbol{a}_{i}\boldsymbol{,Q}_{A_{i}}]\boldsymbol{e}_{\mathbf{m}%
}& =\boldsymbol{a}_{i}\boldsymbol{Q}_{A_{i}}\boldsymbol{e}_{\mathbf{m}}-%
\boldsymbol{Q}_{A_{i}}\boldsymbol{a}_{i}\boldsymbol{e}_{\mathbf{m}}  \notag
\\
& \hspace{-0.5in}={\boldsymbol{a}}_{i}(m_{i}\hbar \,\boldsymbol{e}_{\mathbf{m%
}})-\boldsymbol{Q}_{A_{i}}\boldsymbol{e}%
_{(m_{1},...m_{i-1},m_{i}-1,m_{i+1},...,m_{n})}  \notag \\
& \hspace{-0.5in}=\hbar \,\boldsymbol{e}%
_{(m_{1},...m_{i-1},m_{i}-1,m_{i+1},...,m_{n})}=\hbar \,\boldsymbol{a}_{i}%
\boldsymbol{e}_{\mathbf{m}}\text{.}  \notag
\end{align}%
Hence, 
\begin{equation}
\lbrack \boldsymbol{a}_{i}\boldsymbol{,Q}_{A_{i}}]=\hbar \,\boldsymbol{a}_{i}%
\text{.}  \label{commutation}
\end{equation}%
Moreover, $i\neq j$ implies $[\boldsymbol{a}_{i}\boldsymbol{,Q}_{A_{j}}]=0$.
Taking the adjoint, of the preceding equations, we get 
\begin{equation*}
\lbrack \boldsymbol{a}_{i}^{\dagger }\boldsymbol{,Q}_{A_{i}}]=-\hbar 
\boldsymbol{a}_{i}^{\dagger }\mspace{10mu}\text{and}\mspace{10mu}[%
\boldsymbol{a}_{i}^{\dagger }\boldsymbol{,Q}_{A_{j}}]=0\mspace{10mu}\text{%
whenever}\mspace{10mu}i\neq j\text{. }
\end{equation*}%
If $f_{j}$ is a smooth function on $P$ such that 
\begin{equation}
\{f_{j},A_{k}\}=i\delta _{kj}f_{j},  \label{bracket}
\end{equation}%
then, we can interpret the operator $\boldsymbol{a}_{j}$ as the quantum
operator corresponding to $f_{j}$. In other words, we set $\boldsymbol{a}%
_{j}=\boldsymbol{Q}_{f_{j}}$. This choice is consistent with Dirac's
quantization relations (\ref{Dirac}) 
\begin{equation}
\lbrack \boldsymbol{Q}_{f},\boldsymbol{Q}_{h}]=-i\hbar \,\boldsymbol{Q}%
_{\{f,h\}}
\end{equation}%
because (\ref{commutation}) yields%
\begin{equation*}
\lbrack \boldsymbol{Q}_{f_{j}},\boldsymbol{Q}_{A_{k}}]=-i\hbar \,\boldsymbol{%
Q}_{\{f_{j},A_{k}\}}=-i\hbar \,\boldsymbol{Q}_{(i\delta _{kj}f_{j})}=\delta
_{kj}\hbar \,\boldsymbol{Q}_{f_{j}}\text{.}
\end{equation*}%
Clearly, the function $f_{j}$ is defined by equation (\ref{bracket}) up to
an arbitrary function which commutes with all actions $A_{1},\ldots ,A_{n}$.
Hence, there is a choice involved. We shall use this freedom of choice to
obtain simple expressions for the quantum operators corresponding to the
functions $f_{1},...,f_{n}$. \medskip 

Since $\omega _{\mid U}=\sum_{i=1}^{n}\,\mathop{\! \, \mathrm{d} \!}\nolimits
A_{i}\wedge \mathop{\! \, \mathrm{d} \!}\nolimits\varphi _{i}$, it follows
that the Poisson bracket of $e^{i\varphi _{j}}$ and $A_{k}$ is 
\begin{equation}
\{{\mathrm{e}}^{i\varphi _{j}},A_{k}\}=X_{A_{k}}{\mathrm{e}}^{i\varphi _{j}}=%
\tfrac{\partial }{\partial \varphi _{k}}{\mathrm{e}}^{i\varphi _{j}}=i\delta
_{kj}{\mathrm{e}}^{i\varphi _{j}}.  \label{bracket2}
\end{equation}%
Comparing equations (\ref{bracket}) and (\ref{bracket2}) we see that we may
make the following identification $\boldsymbol{a}_{k}=\boldsymbol{Q}%
_{e^{i\varphi _{k}}}$. Hence, $\boldsymbol{a}_{k}^{\dagger }=\boldsymbol{Q}%
_{e^{-i\varphi _{k}}}$. \medskip 

With this identification we can quantize the following functions on $%
(U,\omega _{\mid U})$.

\begin{itemize}
\item The actions $A_{j}$, $j=1, \ldots n$.

\item The functions ${\mathrm{e}}^{i \varphi _{j}}$, $j =1, \ldots , n$ and
their complex conjugate ${\mathrm{e}}^{-i \varphi _{j}}$, $j=1, \ldots n$,
respectively.

\item Since the operators $\boldsymbol{Q}_{A_{j}}$, $j=1, \ldots , n$,
commute with each other, for any analytic function $H$ of $n$-variables, we
can define an operator $\boldsymbol{Q}_{H(A_{1},...,A_{n})}= H(\boldsymbol{Q}%
_{A_{1}},...,\boldsymbol{Q}_{A_{n}})$ .

\item Since, $\cos \varphi _{j}=%
\mbox{$\frac{\scriptstyle 1}{\scriptstyle
2}\,$}({\mathrm{e}}^{i\varphi _{j}}+{\mathrm{e}}^{-i\varphi _{j}})$ and $%
\sin \varphi _{j}=\mbox{$\frac{\scriptstyle 1}{\scriptstyle 2}\,$}({\mathrm{e%
}}^{i\varphi _{j}}-{\mathrm{e}}^{-i\varphi _{j}})$ we may set $Q_{\cos
\varphi _{j}}=\mbox{$\frac{\scriptstyle 1}{\scriptstyle 2}\,$}(\boldsymbol{a}%
_{j}+\boldsymbol{a}_{j}^{\dagger })$ and $Q_{\sin \varphi _{j}}=%
\mbox{$\frac{{\scriptstyle 1}}{{\scriptstyle 2i }}$}(\boldsymbol{a}_{j}-%
\boldsymbol{a}_{j}^{\dagger })$ .

\item Since the operators $Q_{\cos \varphi _{j}}$ for $j=1, \ldots n$
commute with each other, we can quantize any analytic function of $\cos
\varphi _{j}$ for $j =1, \ldots , n$. Similarly, we can quantize any
analytic function of $\sin \varphi _{j}$ for $j=1, \ldots , n$.

\item We can also quantize functions linear in the actions. For example, 
\begin{eqnarray*}
\boldsymbol{Q}_{A_{i}\cos \varphi _{i}} &=&%
\mbox{$\frac{\scriptstyle
1}{\scriptstyle 2}\,$}(\boldsymbol{Q}_{A_{i}}\boldsymbol{Q}_{\cos \varphi
_{i}}+\boldsymbol{Q}_{\cos \varphi _{i}}\boldsymbol{Q}_{A_{i}}), \\
\rule{0pt}{16pt}\boldsymbol{Q}_{A_{i}\sin \varphi _{i}} &=&%
\mbox{$\frac{\scriptstyle 1}{\scriptstyle 2}\,$}(\boldsymbol{Q}_{A_{i}}%
\boldsymbol{Q}_{\sin \varphi _{i}}+\boldsymbol{Q}_{\sin \varphi _{i}}%
\boldsymbol{Q}_{A_{i}}).
\end{eqnarray*}%
Here the order of the operators on the right hand side is determined by the
requirement that quantization of a real function yields a symmetric operator.
\end{itemize}

\noindent It should be noted that quantization of functions involving the
angles $\varphi _{j}$ for $j=1,\ldots ,n$ gives rise to operators on ${%
\mathfrak{H}}_U$ that are presented as matrices with respect to the basis $%
(e_{\mathbf{m}})$. \medskip

The results described above give a quantization of the symplectic manifold $%
(U, \omega _{\mid U})$ with respect to the real polarization $D_{\mid U}$,
provided $S_U$ is unbounded in every direction. This requirement is
equivalent to the statement that the lattice corresponding to $S_U$ is ${%
\mathbb{Z}}^{n}$. If the lattice $S_U$ is bounded in any direction, we have
to take it into account in our definition of the shifting operators.
Similarly, if the boundary of $U$ contains some Bohr-Sommerfeld tori, we
also have to modify the definition of the shifting operators. These
modifications will be described in the examples treated below.

\section{The $\mathbf{1}$-dimensional harmonic oscillator}


The phase space of a harmonic oscillator is $P={\mathbb{R}}^{2}$ with
coordinates $p$ and $q$ and symplectic form $\omega = 
\mathop{\! \,
\mathrm{d} \!}\nolimits p \wedge \mathop{\! \, \mathrm{d} \!}\nolimits q= %
\mathop{\! \, \mathrm{d} \!}\nolimits \theta $, where $\theta =p\, %
\mathop{\! \, \mathrm{d} \!}\nolimits q $. The Hamiltonian of a harmonic
oscillator is $H= \mbox{$\frac{\scriptstyle 1}{\scriptstyle 2}\,$}
(p^{2}+q^{2})$. The Hamiltonian vector field $X_{H}$ of $H$ is 
\begin{equation*}
X_{H}=\mbox{$\frac{{\scriptstyle \partial H}}{{\scriptstyle \partial p}}$}%
\mbox{$\frac{{\scriptstyle \partial }}{{\scriptstyle \partial q}}$}- %
\mbox{$\frac{{\scriptstyle \partial H}}{{\scriptstyle \partial q}}$}%
\mbox{$\frac{{\scriptstyle \partial }}{{\scriptstyle \partial p}}$}= p%
\mbox{$\frac{{\scriptstyle \partial }}{{\scriptstyle \partial q}}$}-q%
\mbox{$\frac{{\scriptstyle \partial }}{{\scriptstyle \partial p}}$}.
\end{equation*}%
The flow of $X_H$ defines an $\mathop{\mathrm{SU}}\nolimits (1)$-action on $%
P $ given by 
\begin{equation*}
\begin{array}{l}
\Phi :\mathop{\mathrm{SU}}\nolimits (1)\times P\rightarrow P: \\ 
\rule{0pt}{16pt} \hspace{.25in} \big( {\mathrm{e}}^{i \varphi }, (p,q) \big) %
\mapsto \Phi _{{\mathrm{e}}^{i \varphi}}(p,q) = (p\cos \varphi -q\sin
\varphi ,p\sin \varphi +q\cos \varphi ).%
\end{array}%
\end{equation*}
The origin $(0,0)\in {\mathbb{R}}^{2}$ is a fixed point of $\Phi $. The
orbits of $\mathop{\mathrm{SU}}\nolimits (1)$ give rise a Lagrangian
fibration by circles on $U= {\mathbb{R}}^{2}\setminus (0,0)$. \medskip

Using polar coordinates $(r,\varphi )$ on $U$, we have $p=r\cos \varphi$ and 
$q=r\sin \varphi $. Then $\omega =r\, \mathop{\! \, \mathrm{d} \!}\nolimits
r\wedge \mathop{\! \, \mathrm{d} \!}\nolimits \varphi = 
\mathop{\! \,
\mathrm{d} \!}\nolimits \, (\mbox{$\frac{\scriptstyle 1}{\scriptstyle 2}\,$}
r^{2}\, \mathop{\! \, \mathrm{d} \!}\nolimits \varphi )$. This implies that $%
H=\mbox{$\frac{\scriptstyle 1}{\scriptstyle 2}\,$} (p^{2}+q^{2})= %
\mbox{$\frac{\scriptstyle 1}{\scriptstyle 2}\,$} r^{2}$ is an action and $%
\varphi $ is the corresponding angle variable. So polar coordinates $%
(r,\varphi )$ are action-angle variables on $U$. \medskip

The shifting operators $\boldsymbol{a}$ and $\boldsymbol{a}^{\dagger }$,
introduced in \S 4 are 
\begin{equation*}
\boldsymbol{ae}_{m}=\boldsymbol{e}_{m-1}\text{\quad and \quad } \boldsymbol{a%
}^{\dagger }\boldsymbol{e}_{m} =\boldsymbol{e}_{m+1} \text{, \quad for }m > 0
\end{equation*}%
and they correspond to quantum operators $\boldsymbol{a}= \boldsymbol{Q}_{{%
\mathrm{e}}^{i \varphi }}$ and $\boldsymbol{a}^{\dagger }=\boldsymbol{Q}_{{%
\mathrm{e}}^{-i \varphi }}$. The functions ${\mathrm{e}}^{\pm i \varphi }$
do not extend smoothly to the origin $(0,0)$ in ${\mathbb{R}}^{2}$. However,
the functions $z=p+i q=r{\mathrm{e}}^{i \varphi }$ and $\bar{z}=p-i q=r{%
\mathrm{e}}^{-i \varphi }$ are smooth on ${\mathbb{R}}^{2}$. They satisfy
the required Poisson bracket relations%
\begin{eqnarray*}
\{z,H\} &=&%
\mbox{$\frac{{\scriptstyle \partial }}{{\scriptstyle \partial
\varphi }}$}r{\mathrm{e}}^{i \varphi }= i r{\mathrm{e}}^{i \varphi }=i \, z,
\\
\{\bar{z},H\} &=&%
\mbox{$\frac{{\scriptstyle \partial }}{{\scriptstyle
\partial \varphi }}$}r{\mathrm{e}}^{-i \varphi } =-i r{\mathrm{e}}^{-i
\varphi }=-i \, \bar{z}.
\end{eqnarray*}
Therefore, we may introduce new operators $\boldsymbol{b}=\boldsymbol{Q}_{z}$
and $\boldsymbol{b}^{\dagger }=\boldsymbol{Q}_{\bar{z}}$. Equation (\ref%
{Dirac}) yields 
\begin{equation*}
\lbrack \boldsymbol{Q}_{z},\boldsymbol{Q}_{H}] = -i \hbar \, \boldsymbol{Q}%
_{\{z,H\}}=\hbar \, \boldsymbol{Q}_{z} \mspace{10mu}\mathrm{and} %
\mspace{10mu} \lbrack \boldsymbol{Q}_{\bar{z}},\boldsymbol{Q}_{H}] =-i \hbar
\, \boldsymbol{Q}_{\{\bar{z},H\}}=-\hbar \, \boldsymbol{Q}_{\bar{z}} .
\end{equation*}
In other words, $[ \boldsymbol{b},{\boldsymbol{Q}}_{H}] = \hbar \, 
\boldsymbol{b}$ and $[ \boldsymbol{b}^{\dagger },\boldsymbol{Q}_{H}] =
-\hbar \, \boldsymbol{b}^{\dagger }$. Hence, for every $m>0$, we have 
\begin{align}
\boldsymbol{Q}_{H}\boldsymbol{be}_{m} &= \boldsymbol{bQ}_{H}\boldsymbol{e}%
_{m}-\hbar \, \boldsymbol{be}_{m}=(m-1)\hbar \, \boldsymbol{be}_{m},  \notag
\\
\rule{0pt}{16pt} \boldsymbol{Q}_{H}\boldsymbol{b}^{\dagger }\boldsymbol{e}%
_{m} &= \boldsymbol{b}^{\dagger }\boldsymbol{Q}_{H}\boldsymbol{e}_{m}+ \hbar
\, \boldsymbol{b}^{\dagger }\boldsymbol{e}_{m}= (m+1)\hbar \, \boldsymbol{b}%
^{\dagger }\boldsymbol{e}_{m}.  \notag
\end{align}
Observe that $\boldsymbol{b^{\dagger }b}$ commutes with $\boldsymbol{Q}_{H}$%
, because 
\begin{eqnarray*}
\lbrack \boldsymbol{b^{\dagger }b,Q}_{H}] &=&\boldsymbol{b}^{\dagger }%
\boldsymbol{bQ}_{H}-\boldsymbol{Q}_{H}\boldsymbol{b}^{\dagger }\boldsymbol{b}%
=\boldsymbol{b}^{\dagger }\boldsymbol{(bQ}_{H}-\boldsymbol{Q}_{H}\boldsymbol{%
b})\text{ }+\boldsymbol{(b^{\dagger }Q}_{H}-\boldsymbol{Q}_{H}\boldsymbol{b}%
^{\dagger })\boldsymbol{b} \\
&=&\boldsymbol{b}^{\dagger }[\boldsymbol{b},\boldsymbol{Q}_{H}]+[\boldsymbol{%
b^{\dagger },Q}_{H}]\boldsymbol{b}=\boldsymbol{b^{\dagger }(}\hbar \, 
\boldsymbol{b})+(-\hbar \, \boldsymbol{b^{\dagger })b}=0.
\end{eqnarray*}%
In other words, $[\boldsymbol{Q}_{\bar{z}}\boldsymbol{Q}_{z},\boldsymbol{Q}%
_{H}]=0.$ Since $\bar{z}z=r^{2}=2H$, we may assume that 
\begin{equation*}
\boldsymbol{b^{\dagger }b}=\boldsymbol{Q}_{\bar{z}}\boldsymbol{Q}_{z}=%
\boldsymbol{Q}_{\bar{z}z}=2\boldsymbol{Q}_{H}\text{. }
\end{equation*}%
This implies that for every $m\geq 0,$%
\begin{equation*}
\left\Vert \boldsymbol{be}_{m}\right\Vert ^{2}=\langle \boldsymbol{e}%
_{m}\mid \boldsymbol{b^{\dagger }be}_{m}\rangle =\langle \boldsymbol{e}%
_{m}\mid 2\boldsymbol{Q}_{H}\boldsymbol{e}_{m}\rangle =2m\hbar \left\Vert 
\boldsymbol{e}_{m}\right\Vert ^{2}=2m\hbar .
\end{equation*}%
Hence, $\boldsymbol{be}_{0}=0$ and we can choose a normalization factor so
that 
\begin{equation*}
\boldsymbol{be}_{m}=\sqrt{2m\hbar }\, \boldsymbol{e}_{m-1}\text{\quad for }%
~m>0.
\end{equation*}%
Since 
\begin{equation*}
\langle \boldsymbol{e}_{m+1}\mid \boldsymbol{b^{\dagger }e}_{m}\rangle
=\langle \boldsymbol{be}_{m+1}\mid \boldsymbol{e}_{m}\rangle = \sqrt{%
2(m+1)\hbar }\, \langle \boldsymbol{e}_{m}\mid \boldsymbol{e}_{m}\rangle = 
\sqrt{2(m+1)\hbar },
\end{equation*}
we obtain $\boldsymbol{b}^{\dagger }\boldsymbol{e}_{m}=\sqrt{2(m+1)\hbar }\, 
\boldsymbol{e}_{m+1}$. Therefore, for $m>0$ we may write 
\begin{equation}
\boldsymbol{be}_{m}=\sqrt{2m\hbar }\, \boldsymbol{ae}_{m}\text{\quad and
\quad }\boldsymbol{b}^{\dagger }\boldsymbol{e}_{m}=\sqrt{2(m+1)\hbar }\, {%
\boldsymbol{a}}^{\dagger} {\mathbf{e}}_{m}\text{.}  \label{normalization}
\end{equation}
Thus, 
\begin{equation}
\boldsymbol{Q}_{z}e_{m}=\sqrt{2m\hbar }\,\boldsymbol{e}_{m-1}\text{ and }%
\boldsymbol{Q}_{\bar{z}}e_{m}=\sqrt{2(m+1)\hbar }\,\boldsymbol{e}_{m+1}
\label{quantizotion}
\end{equation}
Since $z=p+i q$, we get $p= \mbox{$\frac{\scriptstyle 1}{\scriptstyle 2}\,$}
(z+\bar{z})$ and $q=\mbox{$\frac{{\scriptstyle 1}}{{\scriptstyle 2i}}$}(\bar{%
z}-z)$. Thus the quantization of $p$ and $q$ by operators is $\boldsymbol{Q}%
_{p}=\mbox{$\frac{\scriptstyle 1}{\scriptstyle 2}\,$} (\boldsymbol{b+b}%
^{\dagger })$ and $\boldsymbol{Q}_{q}= 
\mbox{$\frac{\scriptstyle
1}{\scriptstyle 2}\,$} (\boldsymbol{b}^{\dagger }- \boldsymbol{b})$.
Therefore,%
\begin{equation}
\boldsymbol{Q}_{p}\boldsymbol{e}_{m}=\sqrt{%
\mbox{$\frac{{\scriptstyle
m\hbar}}{{\scriptstyle 2}}$}}\, \boldsymbol{e}_{m-1}+\sqrt{%
\mbox{$\frac{{\scriptstyle (m+1)\hbar }}{{\scriptstyle 2}}$}}\, \boldsymbol{e%
}_{m+1}, \quad m > 1; \quad {\mathbf{Q}}_p\, {\mathbf{e}}_1 = \sqrt{%
\mbox{$\frac{{\scriptstyle 2\hbar}}{{\scriptstyle 2}}$}} \, {\mathbf{e}}_2
\label{Qp}
\end{equation}%
and 
\begin{equation}
\boldsymbol{Q}_{q}\boldsymbol{e}_{m}= i\, \sqrt{%
\mbox{$\frac{{\scriptstyle
(m+1)\hbar }}{{\scriptstyle 2}}$}}\boldsymbol{e}_{m+1}-i\, \sqrt{%
\mbox{$\frac{{\scriptstyle m\hbar }}{{\scriptstyle 2}}$}}\, \boldsymbol{e}%
_{m-1}\quad m > 1; \quad {\mathbf{Q}}_q{\mathbf{e}}_1 = i\, \sqrt{%
\mbox{$\frac{{\scriptstyle 2\hbar}}{{\scriptstyle 2}}$}}\, {\mathbf{e}}_2 .
\label{Qq}
\end{equation}%
Equation (\ref{Qp}) can be rewritten in the matrix notation as 
\begin{equation*}
\boldsymbol{Q}_{p}=\left( 
\mbox{\footnotesize $
\begin{array}{ccccc}
0 & \sqrt{\frac{2\hbar }{2}} & 0 & 0 & ... \\ 
\sqrt{\frac{2\hbar }{2}} & 0 & \sqrt{\frac{3\hbar }{2}} & 0 & ... \\ 
0 & \sqrt{\frac{3\hbar }{2}} & 0 & \sqrt{\frac{4\hbar }{2}} & ... \\ 
0 & 0 & \sqrt{\frac{4\hbar }{2}} & 0 & ... \\ 
... & ... & ... & ... & ...\end{array}$} \right)
\end{equation*}%
In a similar way we can write a matrix presentation for $\boldsymbol{Q}_{q}.$

\section{Quantization of coadjoint orbits of $\mathbf{SO(3)}$}


In this section we give the Bohr-Sommerfeld-Heisenberg quantization of
coadjoint orbits of $\mathop{\mathrm{SO}}\nolimits(3)$ on ${%
\mathop{\mathrm{so}}\nolimits(3)}^{\ast }$. Below we show that coadjoint
orbits of $\mathop{\mathrm{SO}}\nolimits(3)$ are spheres in ${\mathbb{R}}%
^{3}.$ Hence, 
\begin{equation*}
P=\{(x^{1},x^{2},x^{3})\in {\mathbb{R}}^{3}\,\mathop{\rule[-4pt]{.5pt}{13pt}%
\, }\nolimits\,(x^{1})^{2}+(x^{2})^{2}+(x^{3})^{2}=r^{2}\}=S_{r}^{2}.
\end{equation*}%
For each $i=1,2,3$, we set \thinspace $J^{i}=x_{\mid P}^{i}$. \medskip 

The following discussion shows that the standard symplectic form\footnote{%
The expression (\ref{SpinOmega}) for the symplectic form on a coadjoint
orbit of $SO(3)$ is the one used by J.-M. Souriau in one of his ledtures.}
on $P$ is%
\begin{equation}
\omega =-\mbox{$\frac{{\scriptstyle 1}}{{\scriptstyle 2r^{2}}}$}%
\sum_{i,j,k=1}^{3}{\varepsilon }_{ijk}J^{i}\,\mathop{\! \, \mathrm{d} \!}%
\nolimits J^{j}\wedge \mathop{\! \, \mathrm{d} \!}\nolimits J^{k}=%
\mbox{$\frac{{\scriptstyle 1}}{{\scriptstyle r}}$}{\mathrm{vol}}_{S_{r}^{2}},
\label{SpinOmega}
\end{equation}%
where ${\mathrm{vol}}_{S_{r}^{2}}$ is the standard volume form on $S_{r}^{2}$
with $\int_{S_{r}^{2}}{\mathrm{vol}}_{S_{r}^{2}}=4\pi r^{2}$. \medskip 

First we recall some basic facts about the Lie algebra $\mathop{\mathrm{so}}%
\nolimits (3)$ of the rotation group $\mathop{\mathrm{SO}}\nolimits (3)$.The
map 
\begin{equation}
j: \mathop{\mathrm{so}}\nolimits (3) \rightarrow {\mathbb{R}}^3:\widehat{X}
= 
\mbox{{\tiny $\begin{pmatrix} 
0 & -x_3 & x_2 \\
x_3 & 0 & -x_1 \\
-x_2 & x_1 & 0  \end{pmatrix} $}} \mapsto x= (x_1,x_2,x_3)
\label{eq-sec6one}
\end{equation}
identifies the Lie algebra $\mathop{\mathrm{so}}\nolimits (3)$ with ${%
\mathbb{R}}^3$. A short calculation shows that $j([\widehat{X}, \widehat{Y}%
]) = x \times y$. Thus $j$ is an isomorphism of the Lie algebra $(%
\mathop{\mathrm{so}}\nolimits (3), [\, , \, ])$ with the Lie algebra $({%
\mathbb{R}}^3, \times)$. It is also an isometry from $(\mathop{\mathrm{so}}%
\nolimits (3), \mathrm{k})$ to $({\mathbb{R}}^3, (\, , \, ) )$, where $%
\mathrm{k}$ is the Killing form on $\mathop{\mathrm{so}}\nolimits (3)$ and $%
(\, , \, ) $ is the Euclidean inner product on ${\mathbb{R}}^3$. To see this
we compute 
\begin{align}
\mathrm{k}(\widehat{X}, \widehat{Y}) & = 
\mbox{$\frac{\scriptstyle
1}{\scriptstyle 2}\,$} \mathrm{tr}\, \widehat{X} {\widehat{Y}}^T = %
\mbox{$\frac{\scriptstyle 1}{\scriptstyle 2}\,$} \mathrm{tr}\, 
\mbox{{\tiny $\begin{pmatrix} 
0 & -x_3 & x_2 \\
x_3 & 0 & -x_1 \\
-x_2 & x_1 & 0  \end{pmatrix} $}}\, 
\mbox{{\tiny $\begin{pmatrix} 
0 & y_3 & -y_2 \\
-y_3 & 0 & y_1 \\
y_2 & y_1 & 0  \end{pmatrix} $}}  \notag \\
& = x_1y_1+x_2y_2+x_3y_3 = (x,y).  \notag
\end{align}
Note that for every $\widehat{X} \in \mathop{\mathrm{so}}\nolimits (3)$ and
every $y \in {\mathbb{R}}^3$ we have 
\begin{equation}
\widehat{X}y = 
\mbox{{\tiny $\begin{pmatrix} 
0 & -x_3 & x_2 \\
x_3 & 0 & -x_1 \\
-x_2 & x_1 & 0  \end{pmatrix} $}} \, 
\mbox{{\tiny $\begin{pmatrix} 
y_1 \\ y_2 \\ y_3 \end{pmatrix} $}} = x \times y.  \label{eq-sec6two}
\end{equation}
For every $\widehat{X}$, $\widehat{Y} \in \mathop{\mathrm{so}}\nolimits (3)$
using (\ref{eq-sec6two}) we can rewrite $j([\widehat{X}, \widehat{Y}]) =
x\times y$ as $j({\mathop{\mathrm{ad}}\nolimits }_{\widehat{X}}\widehat{Y})
= \widehat{X}j(\widehat{Y})$, which is equivalent to $\big( j({%
\mathop{\mathrm{ad}}\nolimits}_{\widehat{X}})j^{-1} \big) j(\widehat{Y}) = 
\widehat{X}j(\widehat{Y})$, that is, 
\begin{equation}
j({\mathop{\mathrm{ad}}\nolimits}_{\widehat{X}})j^{-1} = \widehat{X}, \text{%
\quad for every}\mspace{5mu} \widehat{X} \in \mathop{\mathrm{so}}\nolimits
(3).  \label{eq-sec6twostarstar}
\end{equation}

Let $R \in \mathop{\mathrm{SO}}\nolimits (3)$ and $\widehat{Y} \in %
\mathop{\mathrm{so}}\nolimits (3)$. Then we have 
\begin{equation}
j({\mathop{\mathrm{Ad}}\nolimits}_R \widehat{Y}) = Ry.  \label{eq-sec6three}
\end{equation}
To prove (\ref{eq-sec6three}) we need the following formula, which holds for
any linear Lie group $\mathfrak{G}$ and its associated linear Lie algebra $%
\mathfrak{g}$, namely 
\begin{equation}
{\mathop{\mathrm{Ad}}\nolimits}_{\exp tX} = \exp t\, {\mathop{\mathrm{ad}}%
\nolimits}_X,  \label{eq-sec6four}
\end{equation}
for every $X \in \mathfrak{g}$ and every $t \in \mathbb{R} $. To verify that
(\ref{eq-sec6four}) holds, we note that the right and left hand sides of (%
\ref{eq-sec6four}) are each $1$-parameter subgroups of $\mathfrak{G}$ with
the same tangent vector at $t=0$, namely, ${\mathop{\mathrm{ad}}\nolimits }%
_X $. Therefore the $1$-parameter subgroups are equal. \medskip

Returning to the proof of (\ref{eq-sec6three}), using (\ref{eq-sec6four}) we
get 
\begin{equation*}
j({\mathop{\mathrm{Ad}}\nolimits}_{\exp t\widehat{X}})j^{-1} = \exp (t \, j({%
\mathop{\mathrm{ad}}\nolimits}_{\widehat{X}})j^{-1}) = \exp t\widehat{X},
\end{equation*}
that is, for every $\widehat{Y} \in \mathop{\mathrm{so}}\nolimits (3)$ we
have 
\begin{equation}
j\big( ({\mathop{\mathrm{Ad}}\nolimits}_{\exp t \widehat{X}})\widehat{Y}) =
( \exp t \widehat{X}) j(\widehat{Y}).  \label{eq-sec6six}
\end{equation}
Since $\mathop{\mathrm{SO}}\nolimits (3)$ is compact and connected, for
every $R \in \mathop{\mathrm{SO}}\nolimits (3)$ there is a $\widehat{X} \in %
\mathop{\mathrm{so}}\nolimits (3)$ such that $R = \exp \widehat{X}$. Thus (%
\ref{eq-sec6six}) implies that for every $R \in SO (3)$ equation (\ref%
{eq-sec6three}) holds. So (\ref{eq-sec6three}) is an integrated version of (%
\ref{eq-sec6twostarstar}). \medskip

Now we calculate the standard symplectic form on an $\mathop{\mathrm{SO}}%
\nolimits (3)$-adjoint orbit. The $\mathop{\mathrm{SO}}\nolimits (3)$%
-adjoint orbit through $\widehat{J} \in \mathop{\mathrm{so}}\nolimits (3)$
is ${\mathcal{O}}_{\widehat{J}} = \{ {\mathop{\mathrm{Ad}}\nolimits}_R%
\widehat{J} \in \mathop{\mathrm{so}}\nolimits (3) \, \, \mathop{%
\rule[-4pt]{.5pt}{13pt}\, }\nolimits \, R \in \mathop{\mathrm{SO}}\nolimits
(3) \}$. The standard symplectic form $\Omega $ on ${\mathcal{O}}_{\widehat{J%
}}$ is 
\begin{equation}
\Omega (\widehat{J})(X^{\widehat{\xi}}(\widehat{J}), X^{\widehat{\zeta }}(%
\widehat{J})) = - \mathrm{k}(\widehat{J}, [ \widehat{\xi}, \widehat{\zeta }
]),  \label{eq-sec6s2one}
\end{equation}
where $\widehat{\xi}, \widehat{\zeta } \in \mathop{\mathrm{so}}\nolimits (3)$
and $X^{\widehat{\eta }}(\widehat{J}) = -{\mathop{\mathrm{ad}}\nolimits}_{%
\widehat{J}}\widehat{\eta } = -[\widehat{J}, \widehat{\eta }]$, which
defines a vector field on ${\mathcal{O}}_{\widehat{J}}$ for every $\widehat{%
\eta } \in \mathop{\mathrm{so}}\nolimits (3)$. Because 
\begin{equation*}
{\mathop{\mathrm{Ad}}\nolimits}_R X^{\widehat{\eta }}(\widehat{J}) = - {%
\mathop{\mathrm{Ad}}\nolimits}_R[\widehat{J}, \widehat{\eta }] = - [{%
\mathop{\mathrm{Ad}}\nolimits}_R\widehat{J}, {\mathop{\mathrm{Ad}}\nolimits}%
_R\widehat{\eta }] = X^{{\mathop{\mathrm{Ad}}\nolimits}_R\widehat{\eta }}({%
\mathop{\mathrm{Ad}}\nolimits}_R \widehat{J}),
\end{equation*}
we get 
\begin{align}
\Omega ({\mathop{\mathrm{Ad}}\nolimits}_R\widehat{J})({\mathop{\mathrm{Ad}}%
\nolimits}_RX^{\widehat{\xi }}(\widehat{J}), {\mathop{\mathrm{Ad}}\nolimits}%
_RX^{\widehat{\zeta }}(\widehat{J})) & = \Omega ({\mathop{\mathrm{Ad}}%
\nolimits}_R\widehat{J})( X^{{\mathop{\mathrm{Ad}}\nolimits}_R\widehat{\xi }%
}({\mathop{\mathrm{Ad}}\nolimits}_R \widehat{J}), X^{{\{Ad}_R\widehat{\zeta }%
} ({\mathop{\mathrm{Ad}}\nolimits}_R \widehat{J}))  \notag \\
& \hspace{-1in} = -\mathrm{k}({\mathop{\mathrm{Ad}}\nolimits}_R\widehat{J},
[ {\mathop{\mathrm{Ad}}\nolimits}_R\widehat{\xi}, {\mathop{\mathrm{Ad}}%
\nolimits}_R\widehat{\zeta } ]) = -\mathrm{k}({\mathop{\mathrm{Ad}}\nolimits}%
_R\widehat{J}, {\mathop{\mathrm{Ad}}\nolimits}_R[\widehat{\xi }, \widehat{%
\zeta }])  \notag \\
& \hspace{-1in} = -\mathrm{k}(\widehat{J}, [\widehat{\xi}, \widehat{\zeta }
]) = \Omega (\widehat{J})(X^{\widehat{\xi}}(\widehat{J}), X^{\widehat{\zeta }%
} (\widehat{J})) .  \notag
\end{align}
This shows that $\Omega$ is a $2$-form on ${\mathcal{O}}_{\widehat{J}}$. It
is closed since ${\mathcal{O}}_{\widehat{J}}$ is a $2$-dimensional smooth
manifold. It is nondegenerate for if $0 = \Omega (\widehat{J})(X^{\widehat{%
\xi}}(\widehat{J}), X^{\widehat{\zeta }}(\widehat{J}))$ for every $X^{%
\widehat{\zeta }}$ with $\widehat{\zeta } \in \mathop{\mathrm{so}}\nolimits
(3)$, then we obtain $0 = \mathrm{k}(\widehat{J}, [ \widehat{\xi}, \widehat{%
\zeta } ] ) = \mathrm{k}([\widehat{J}, \widehat{\xi}], \widehat{\zeta})$ for
every $\widehat{\zeta } \in \mathop{\mathrm{so}}\nolimits (3)$. Since $%
\mathrm{k}$ is nondegenerate, this implies that $[\widehat{J}, \widehat{\xi}
] =0$. But then $X^{\widehat{\xi}}(\widehat{J}) = - [\widehat{J}, \widehat{%
\xi}] =0$. \medskip

Using the bijection $j$ (\ref{eq-sec6one}) to identify ${\mathop{\mathrm{Ad}}%
\nolimits}_R\widehat{J}$ with $R\widehat{J}$ by (\ref{eq-sec6three}), we see
that the $\mathop{\mathrm{SO}}\nolimits (3)$-adjoint orbit ${\mathcal{O}}_{%
\widehat{J}}$ may be identified with the $2$-sphere $S^2_r = \{ R J \in {{%
\mathbb{R}} }^3 \, \, \mathop{\rule[-4pt]{.5pt}{13pt}\, }\nolimits \, R \in %
\mathop{\mathrm{SO}}\nolimits (3) \}$. Here $r^2 = (J,J)$. We may rewrite
the definition of $\Omega $ (\ref{eq-sec6s2one}) as 
\begin{equation*}
\Omega (\widehat{J})(-[\widehat{J}, \widehat{\xi}], -[\widehat{J}, \widehat{%
\zeta}]) = - \mathrm{k}(\widehat{J}, [\widehat{\xi}, \widehat{\zeta }]).
\end{equation*}
Thus we may identify $\Omega $ with the symplectic form $\omega $ on $S^2_r$
given by 
\begin{equation}
\omega (J)(-J \times \xi , -J\times \zeta ) = -(J, \xi \times \zeta ).
\label{eq-sec6s2two}
\end{equation}
Note that $J\times \xi $ and $J \times \zeta $ both lie in $T_JS^2_r$.
\medskip

The vector field $X^{\widehat{\eta }}$ on ${\mathcal{O}}_{\widehat{J}}$,
defined by $X^{\widehat{\eta }}(\widehat{J}) = -[\widehat{J}, \widehat{\eta }%
] $, corresponds to the vector field $X^{\eta }$ on $S^2_r$ defined by $%
X^{\eta }(J) = - J \times \eta $, because the curve $t \mapsto {%
\mathop{\mathrm{Ad}}\nolimits}_{\exp t\widehat{\eta }}\widehat{J}$ in ${%
\mathcal{O}}_{\widehat{J}}$ is identified under the map $j$ with the curve $%
t \mapsto (\exp t\eta )J$ on $S^2_r$. Therefore the tangent vector $X^{%
\widehat{\eta }}(\widehat{J})$ at $\widehat{J}$ corresponds to the tangent
vector $X^{\eta }(J)$ at $J$, namely $\widehat{\eta }(J) = -J \times \eta $.
So we may rewrite the definition of $\omega $ (\ref{eq-sec6s2two}) as 
\begin{equation}
\omega (J)(X^{\xi}(J), X^{\zeta }(J)) = - (J, \xi \times \zeta ).
\label{eq-sec6s2three}
\end{equation}

Next we show that (\ref{SpinOmega}) holds. Evaluating the left hand side of (%
\ref{SpinOmega}) on the tangent vectors $X^{\xi}(J)$ and $X^{\zeta}(J)$
gives 
\begin{align}
-\mbox{$\frac{{\scriptstyle 1}}{{\scriptstyle r^2}}$} \sum^3_{i=1} J^i \, %
\mbox{$\frac{\scriptstyle 1}{\scriptstyle 2}\,$} \sum^3_{j,k =1} {%
\varepsilon }_{ijk} (\mathop{\! \, \mathrm{d} \!}\nolimits J^j \wedge %
\mathop{\! \, \mathrm{d} \!}\nolimits J^k)(X^{\xi}(J), X^{\zeta }(J)) & = 
\notag \\
& \hspace{-3in} = -\mbox{$\frac{{\scriptstyle 1}}{{\scriptstyle r^2}}$}
\sum^3_{i,j,k =1} J^i \, {\varepsilon }_{ijk} \mathop{\! \, \mathrm{d} \!}%
\nolimits J^j(X^{\xi}(J))\, \mathop{\! \, \mathrm{d} \!}\nolimits
J^k(X^{\zeta }(J))  \notag \\
& \hspace{-3in} = -\mbox{$\frac{{\scriptstyle 1}}{{\scriptstyle r^2}}$} ( J,
(\xi \times J) \times (\zeta \times J)), \quad 
\mbox{see
(\ref{eq-sec6s2five}) below}  \notag \\
&\hspace{-3in} = \mbox{$\frac{{\scriptstyle 1}}{{\scriptstyle r^2}}$}
(J,J)(\xi \times J, \zeta ) = -(J, \xi \times \zeta ) = \omega
(J)(X^{\xi}(J), X^{\zeta }(J)).  \notag
\end{align}
This proves (\ref{SpinOmega}) provided that we show 
\begin{equation}
\mathop{\! \, \mathrm{d} \!}\nolimits J^{\ell }(X^{\eta }(J)) = (\eta \times
J)^{\ell }.  \label{eq-sec6s2five}
\end{equation}
By definition $J^{\ell } = x^{\ell }|S^2_r$, where $x^{\ell }$ is the ${\ell 
}^{\mathrm{th}}$ coordinate function on ${\mathbb{R} }^3$. Now 
\begin{align}
\mathop{\! \, \mathrm{d} \!}\nolimits J^{\ell }(X^{\eta }(J)) & = 
\mbox{${\displaystyle \frac{\dee}{\dee t}}
\rule[-10pt]{.5pt}{25pt} \raisebox{-10pt}{$\, {\scriptstyle t=0}$}$} J^{\ell}%
\big( (\exp t\widehat{\eta })J \big)  \notag \\
& = J^{\ell }\big( 
\mbox{${\displaystyle \frac{\dee}{\dee t}}
\rule[-10pt]{.5pt}{25pt} \raisebox{-10pt}{$\, {\scriptstyle t=0}$}$} (\exp t%
\widehat{\eta })J \big) , \quad 
\mbox{since $J^{\ell }$ is 
a linear function on ${\R }^3$}  \notag \\
& = J^{\ell }(\widehat{\eta }J) =J^{\ell }(\eta \times J) = (\eta \times
J)^{\ell }.  \notag
\end{align}
This completes the verification of (\ref{eq-sec6s2five}) and thus the proof
of (\ref{SpinOmega}). \medskip

Our aim is to obtain an irreducible unitary representation of $%
\mathop{\mathrm{SO}}\nolimits (3)$ corresponding to quantizable coadjoint
orbit. We shall do it in the framework of geometric quantization as
described in \S 2. First we obtain quantum operators ${\mathbf{Q}}_{J^{1}}$, 
${\mathbf{Q}}_{J^{2}}$, ${\mathbf{Q}}_{J^{3}}$. Next, we show that the
rescaled operators $\mbox{$\frac{{\scriptstyle i }}{{\scriptstyle \hbar }}$}%
\, {\mathbf{Q}}_{J^{1}}$, $%
\mbox{$\frac{{\scriptstyle i }}{{\scriptstyle
\hbar }}$}\, {\mathbf{Q}}_{J^{2}}$ and $%
\mbox{$\frac{{\scriptstyle i
}}{{\scriptstyle \hbar }}$}\, {\mathbf{Q}}_{J^{3}}$ give rise to a
representation of $\mathop{\mathrm{so}}\nolimits (3)$. We could proceed
directly by setting $\hbar =1$, or even $\hbar =i $, but in this way we
would lose the connection between geometric quantization in mechanics and in
representation theory. \medskip

We assume that $(P,\omega )$ is prequantizable. This means that $%
\int_{P}\omega =nh$, where $n\in \mathbb{Z}$. Introducing spherical polar
coordinates 
\begin{equation*}
J^{1}=r\sin \theta \cos \varphi ,~J^{2}=r\sin \theta \sin \varphi
,~J^{3}=r\cos \theta
\end{equation*}%
on $S_{r}^{2}$ we get $\omega =r\sin \theta \mathop{\! \, \mathrm{d} \!}%
\nolimits\varphi \wedge \mathop{\! \, \mathrm{d} \!}\nolimits\theta =%
\mbox{$\frac{{\scriptstyle 1}}{{\scriptstyle r}}$}{\mathrm{vol}}_{S_{r}^{2}}$%
. Hence, 
\begin{equation*}
\int_{P}\omega =r\int_{0}^{2\pi }\mathop{\! \, \mathrm{d} \!}%
\nolimits\varphi \int_{0}^{\pi }\sin \theta \,\mathop{\! \, \mathrm{d} \!}%
\nolimits\theta =4\pi r,
\end{equation*}%
and the integrality condition reads $4\pi r=nh$. Equivalently, $r=%
\mbox{$\frac{{\scriptstyle n}}{{\scriptstyle 2}}$}\hbar $ where $\hbar =%
\mbox{$\frac{{\scriptstyle h}}{{\scriptstyle 2\pi}}$}$. Next 
\begin{equation*}
X_{J^{3}}\mbox{$\, \rule{8pt}{.5pt}\rule{.5pt}{6pt}\, \, $}r\sin \theta \,%
\mathop{\! \, \mathrm{d} \!}\nolimits\varphi \wedge \mathop{\! \, \mathrm{d}
\!}\nolimits\theta =r\sin \theta \,\mathop{\! \, \mathrm{d} \!}%
\nolimits\theta =-\mathop{\! \, \mathrm{d} \!}\nolimits J^{3}
\end{equation*}%
implies that $X_{J^{3}}=%
\mbox{$\frac{{\scriptstyle \partial }}{{\scriptstyle
\partial \varphi }}$}$. Thus, the integral curves of $X_{J^{3}}$ are circles 
$J^{3}=\mathrm{const.}$ They define the leaves of a singular real
polarization of $S_{{r}^{2}}$ with singularities poles at $J^{3}=\pm r=\pm %
\mbox{$\frac{{\scriptstyle n}}{{\scriptstyle 2}}$}\hbar $. Locally, we have 
\begin{equation}
\omega =r\sin \theta \,\mathop{\! \, \mathrm{d} \!}\nolimits\varphi \wedge %
\mathop{\! \, \mathrm{d} \!}\nolimits\theta =\mathop{\! \, \mathrm{d} \!}%
\nolimits\,(r\cos \theta \,\mathop{\! \, \mathrm{d} \!}\nolimits\varphi )=%
\mathop{\! \, \mathrm{d} \!}\nolimits\,(J^{3}\,\mathop{\! \, \mathrm{d} \!}%
\nolimits\varphi ).  \label{action-angle}
\end{equation}%
Thus $(J^{3},\varphi )$ are action-angle coordinates for our integrable
system $(J^{3},S_{{r}^{2}},\omega )$. Since $r=%
\mbox{$\frac{{\scriptstyle
n}}{{\scriptstyle 2}}$}\hbar $, the Bohr-Sommerfeld conditions 
\begin{equation}
\int_{J^{3}=\mathrm{const.}}r\cos \theta \,\mathop{\! \, \mathrm{d} \!}%
\nolimits\varphi =mh,  \label{2}
\end{equation}%
read 
\begin{equation*}
\int_{0}^{2\pi }\mbox{$\frac{{\scriptstyle n}}{{\scriptstyle 2}}$}\hbar \cos
\theta \mathop{\! \, \mathrm{d} \!}\nolimits\varphi =2\pi (%
\mbox{$\frac{{\scriptstyle n}}{{\scriptstyle 2}}$}\hbar \cos \theta )=mh%
\text{,}
\end{equation*}%
which implies that $\mbox{$\frac{{\scriptstyle n}}{{\scriptstyle 2}}$}\hbar
\cos \theta =m\hbar $ or $\cos \theta =2%
\mbox{$\frac{{\scriptstyle
m}}{{\scriptstyle n}}$}$. Since $-1\leq \cos \theta \leq 1$, it follows that 
$-1\leq 2\mbox{$\frac{{\scriptstyle M}}{{\scriptstyle n}}$}\leq 1$ or $-%
\mbox{$\frac{{\scriptstyle n}}{{\scriptstyle 2}}$}\leq m\leq %
\mbox{$\frac{{\scriptstyle n}}{{\scriptstyle 2}}$}$. We now assume that $s=%
\mbox{$\frac{{\scriptstyle n}}{{\scriptstyle 2}}$}$ is an integer. Then $%
-s\leq m\leq s$. Thus we get a family $\theta _{m}$ of angles in spherical
coordinates on $S_{{s\hbar }^{2}}$ for which 
\begin{equation}
\cos {\theta }_{m}=\mbox{$\frac{{\scriptstyle m}}{{\scriptstyle s}}$}\text{%
,\quad where}\,-s\leq m\leq s.  \label{0}
\end{equation}%
For $m=\pm s$, we get the north pole $(0,0,s\hbar )$ and the south pole $%
(0,0,-s\hbar )$ of $S_{{s\hbar }^{2}}$. These are the singular points of our
Bohr-Sommerfeld set.\medskip

Let $(\boldsymbol{e}_{m})$ be a basis of $\mathfrak{H}$ consisting of
eigenvectors of $\boldsymbol{Q}_{J_{3}}$. For each integer $m$ between $-s$
and $s$, we have 
\begin{equation*}
\boldsymbol{Q}_{J_{3}}\boldsymbol{e}_{m}=r\cos \theta _{m}\boldsymbol{e}%
_{m}=s\hbar \cos \theta _{m}\boldsymbol{e}_{m}=s\hbar \mbox{$\frac{{%
\scriptstyle m}}{{\scriptstyle s}}$} \, \boldsymbol{e}_{m}=m\hbar \, 
\boldsymbol{e}_{m},
\end{equation*}%
using equation (\ref{0}). We assume that 
\begin{equation}
(\boldsymbol{e}_{m^{\prime }}\mid \boldsymbol{e}_{m})=\delta _{m^{\prime },
m}.  \label{product}
\end{equation}

Note that the Bohr-Sommerfeld conditions not only give the directions of the
basis vectors $\boldsymbol{e}_{m}$ in $\mathfrak{H}$, but also their
ordering $m\mapsto \boldsymbol{e}_{m}$. As in \S 4, we can define the
shifting operators $\boldsymbol{a}$ and $\boldsymbol{a}^{\dagger }$ on $%
\mathfrak{H}$ by 
\begin{equation}
\boldsymbol{a\, e}_{m}=\boldsymbol{e}_{m-1}\text{\quad and \quad }%
\boldsymbol{a}^{\dagger }\boldsymbol{e}_{m}=\boldsymbol{e}_{m+1}\text{.}
\label{5}
\end{equation}%
As before, we can make an identification $\boldsymbol{a}_{i}=\boldsymbol{Q}_{%
{\mathrm{e}}^{i \varphi }}$ and $\boldsymbol{a}_{i}^{\dagger }= \boldsymbol{Q%
}_{{\mathrm{e}}^{-i \varphi }}$. The functions ${\mathrm{e}}^{i \varphi }$
and ${\mathrm{e}}^{-i \varphi }$ do not extend to the singular points $%
(0,0,-s\hbar )$ and $(0,0,s\hbar )$ of the polarization, which correspond to 
$m=-s$ and $m=s$, respectively. However, the function 
\begin{align}
J_{-}& =\sqrt{r^{2}-(J^{3})^{2}}\,{\mathrm{e}}^{i \varphi }=\sqrt{%
r^{2}-(J^{3})^{2}}\,\cos \varphi +i \sqrt{r^{2}-(J^{3})^{2}}\, \sin \varphi 
\notag \\
& =J^{1}+i J^{2}  \notag
\end{align}%
extends smoothly to the singular points. Similarly, the function 
\begin{equation*}
J_{+}=\sqrt{r^{2}-(J^{3})^{2}}\,{\mathrm{e}}^{-i \varphi }= J^{1}-i J^{2}
\end{equation*}%
extends smoothly to the whole of $S_{{s\hbar }^{2}}$. Moreover, we have 
\begin{equation*}
\{J_{-}, J^{3}\}=i J_{+}\text{\quad and \quad }\{J_{+},J^{3}\}= -i J_{-}.
\end{equation*}%
Hence, we can consider our shifting operators to be quantizations of $J_{+}$
and $J_{-}$. In order to define the operators $\boldsymbol{Q}_{J_{+}}$ and $%
\boldsymbol{Q}_{J_{-}}$ we set 
\begin{equation*}
\boldsymbol{Q}_{J_{-}}\boldsymbol{e}_{m}=a_{m}\boldsymbol{e}_{m-1} \text{%
\quad and \quad } \boldsymbol{Q}_{J_{+}}\boldsymbol{e}_{m}=\boldsymbol{Q}%
_{J_{-}}^{\dagger }\boldsymbol{e}_{m}=a_{m+1}\boldsymbol{e}_{m+1},
\end{equation*}%
where the real coefficients $a_{m}$ are to be defined so that $a_{-s}=0$ and 
$a_{s+1}=0$. We have 
\begin{equation*}
\boldsymbol{Q}_{J_{+}}\boldsymbol{Q}_{J_{-}}\boldsymbol{e}_{m}=a_m%
\boldsymbol{Q}_{J_{+}}\boldsymbol{e}_{m+1}=a_m^{2}\boldsymbol{e}_{m}
\end{equation*}%
and 
\begin{equation*}
\boldsymbol{Q}_{J_{-}}\boldsymbol{Q}_{J+}\boldsymbol{e}_{m}=a_{m+1}%
\boldsymbol{Q}_{J_{-}}\boldsymbol{e}_{m-1}=a_{m+1}^{2}\boldsymbol{e}_{m}%
\text{. }
\end{equation*}%
Hence $[\boldsymbol{Q}_{J_{+}},\boldsymbol{Q}_{J_{-}}]\boldsymbol{e}%
_{m}=(a_{m}^{2}-a_{m+1}^{2})\boldsymbol{e}_{m}$. Since 
\begin{equation*}
\{J_{+},J_{-}\}=\{J^{1}-iJ^{2},J^{1}+iJ^{2}\}=2i\{J^{1},J^{2}\}=2i\,J^{3},
\end{equation*}%
it follows that we should have 
\begin{equation}
\lbrack \boldsymbol{Q}_{J_{+}},\boldsymbol{Q}_{J_{-}}]=-i\hbar \,\boldsymbol{%
Q}_{2iJ^{3}}=2\hbar \,\boldsymbol{Q}_{J^{3}}.  \label{eq-coadnew}
\end{equation}%
Therefore, 
\begin{equation}
(a_{m}^{2}-a_{m+1}^{2})\boldsymbol{e}_{m}=2\hbar \,\boldsymbol{Q}_{J^{3}}%
\boldsymbol{e}_{m}=2m\hbar ^{2}\,\boldsymbol{e}_{m}  \label{eq-coeffnew}
\end{equation}%
for every $m=-s,....,s$. Hence, 
\begin{equation*}
a_{m+1}^{2}-a_{m}^{2}=-2m\hbar ^{2}\text{\quad or \quad }%
a_{m}^{2}=a_{m+1}^{2}+2m\hbar ^{2}
\end{equation*}%
For $m=0,$ we have $a_{1}^{2}=a_{0}^{2}$. For $m\geq 1$, we get 
\begin{equation*}
a_{m}^{2}=a_{1}^{2}-2\hbar ^{2}\sum_{k=1}^{m-1}k=a_{1}^{2}-\hbar ^{2}(m-1)m%
\text{,}
\end{equation*}%
and 
\begin{equation*}
a_{-m}^{2}=a_{0}^{2}+2\hbar ^{2}\sum_{k=1}^{m}(-k)=a_{0}^{2}-\hbar ^{2}m(m+1)%
\text{. }
\end{equation*}%
The conditions $a_{s+1}=0$ and $a_{-s}=0$ yield $a_{1}^{2}-\hbar
^{2}s(s+1)=0 $ and $a_{0}^{2}-\hbar ^{2}s(s+1)=0$. Hence, $%
a_{1}^{2}=a_{0}^{2}=\hbar ^{2}s(s+1)$ and 
\begin{eqnarray*}
a_{m}^{2} &=&\hbar ^{2}s(s+1)-\hbar ^{2}(m-1)m\text{,} \\
\rule{0pt}{16pt}a_{-m}^{2} &=&\hbar ^{2}s(s+1)-\hbar ^{2}m(m+1)\text{, }
\end{eqnarray*}%
for $m\geq 1$. Thus, for $m=-1,\ldots ,-s,$ we get 
\begin{equation*}
a_{m}^{2}=\hbar ^{2}s(s+1)+\hbar ^{2}m(-m+1)=\hbar ^{2}s(s+1)-\hbar
^{2}m(m-1).
\end{equation*}%
Therefore, 
\begin{equation*}
a_{m}^{2}=\hbar ^{2}s(s+1)+\hbar ^{2}m(-m+1),\text{\quad for all $%
m=-s,\ldots ,s$}
\end{equation*}%
and 
\begin{eqnarray*}
\boldsymbol{Q}_{J_{-}}\boldsymbol{e}_{m} &=&\hbar \sqrt{s(s+1)-(m-1)m}\,%
\boldsymbol{e}_{m-1} = a_m{\mathbf{e}}_{m-1}, \\
\rule{0pt}{16pt}\boldsymbol{Q}_{J_{+}}\boldsymbol{e}_{m} &=&\boldsymbol{Q}%
_{J_{-}}^{\dagger }\,\boldsymbol{e}_{m}=\hbar \sqrt{s(s+1)-m(m+1)}\,%
\boldsymbol{e}_{m+1} = a_{m+1}{\mathbf{e}}_{m+1}.
\end{eqnarray*}%
So 
\begin{align}
[{\mathbf{Q}}_{J_{+}}, {\mathbf{Q}}_{J_{-}}]{\mathbf{e}}_{m} & = {\mathbf{Q}}%
_{J_{+}}{\mathbf{Q}}_{J_{-}}{\mathbf{e}}_{m}- {\mathbf{Q}}_{J_{-}}{\mathbf{Q}%
}_{J_{+}}{\mathbf{e}}_{m}  \notag \\
& \hspace{-.75in} = \hbar \sqrt{s(s+1)-m(m-1)}{\mathbf{Q}}_{J_{+}}%
\boldsymbol{e}_{m+1} - \hbar \sqrt{s(s+1)-(m+1)m}\, {\mathbf{Q}}_{J_{-}} 
\boldsymbol{e}_{m-1}  \notag \\
& \hspace{-.75in}= {\hbar}^2 \big( s(s+1)-m(m-1)\big) \boldsymbol{e}_m -{%
\hbar }^2\big( s(s+1)-(m+1)m\big) \boldsymbol{e}_m  \notag \\
&\hspace{-.75in} = 2m\, {\hbar}^2\, {\mathbf{e}}_m = 2\hbar \, \boldsymbol{Q}%
_{J^{3}} {\mathbf{e}}_m,  \notag
\end{align}
which verfies that (\ref{eq-coadnew}) holds. Since $J^{1}=%
\mbox{$\frac{{\scriptstyle 1}}{{\scriptstyle 2}}$} (J_{+}-J_{-}) $ and $%
J^{2}=\mbox{$\frac{{\scriptstyle 1}}{{\scriptstyle 2i}}$} (J_{+}+J_{-})$, we
get 
\begin{align}
\boldsymbol{Q}_{J^{1}}& =\mbox{$\frac{{\scriptstyle 1}}{{\scriptstyle 2}}$}
\, {\mathbf{Q}}_{J^{+}}{\mathbf{e}}_m + 
\mbox{$\frac{{\scriptstyle
1}}{{\scriptstyle 2}}$} \, {\mathbf{Q}}_{J^{-}}{\mathbf{e}}_m = %
\mbox{$\frac{\scriptstyle 1}{\scriptstyle 2}\,$} a_{m+1}{\mathbf{e}}_{m+1} + %
\mbox{$\frac{\scriptstyle 1}{\scriptstyle 2}\,$} a_m {\mathbf{e}}_{m-1} 
\notag \\
\boldsymbol{Q}_{J^{2}}& =\mbox{$\frac{{\scriptstyle 1}}{{\scriptstyle 2i}}$}
\, {\mathbf{Q}}_{J^{+}}{\mathbf{e}}_m - 
\mbox{$\frac{{\scriptstyle
1}}{{\scriptstyle 2i}}$} \, {\mathbf{Q}}_{J^{-}}{\mathbf{e}}_m = %
\mbox{$\frac{{\scriptstyle 1}}{{\scriptstyle 2i}}$} a_{m+1}{\mathbf{e}}%
_{m+1} - \mbox{$\frac{{\scriptstyle 1}}{{\scriptstyle 2i}}$} a_m {\mathbf{e}}%
_{m-1}  \notag
\end{align}
The operators $\boldsymbol{Q}_{J^{1}}$, $\boldsymbol{Q}_{J^{2}}$, and $%
\boldsymbol{Q}_{J^{3}}$ satisfy the required commutation relations, namely, 
\begin{align}
[{\mathbf{Q}}_{J^1}, {\mathbf{Q}}_{J^2}] {\mathbf{e}}_m & = {\mathbf{Q}}%
_{J^1}{\mathbf{Q}}_{J^2}{\mathbf{e}}_m - {\mathbf{Q}}_{J^2}{\mathbf{Q}}_{J^1}%
{\mathbf{e}}_m  \notag \\
&\hspace{-.75in} = {\mathbf{Q}}_{J^1}\big( 
\mbox{$\frac{{\scriptstyle
1}}{{\scriptstyle 2i}}$} a_{m+1} {\mathbf{e}}_{m+1} -\mbox{$\frac{{%
\scriptstyle 1}}{{\scriptstyle 2i}}$}a_m {\mathbf{e}}_{m-1} \big) -{\mathbf{Q%
}}_{J^2}\big( \mbox{$\frac{{\scriptstyle 1}}{{\scriptstyle 2}}$} a_{m+1} {%
\mathbf{e}}_{m+1} +\mbox{$\frac{{\scriptstyle 1}}{{\scriptstyle 2}}$}a_m {%
\mathbf{e}}_{m-1} \big)  \notag \\
&\hspace{-.75in} = \mbox{$\frac{{\scriptstyle 1}}{{\scriptstyle 2i}}$}a_{m+1}%
{\mathbf{Q}}_{J^1}{\mathbf{e}}_{m+1} -%
\mbox{$\frac{{\scriptstyle
1}}{{\scriptstyle 2i}}$}a_m{\mathbf{Q}}_{J^1}{\mathbf{e}}_{m-1} -%
\mbox{$\frac{{\scriptstyle 1}}{{\scriptstyle 2}}$}a_{m+1}{\mathbf{Q}}_{J^2}{%
\mathbf{e}}_{m+1} -\mbox{$\frac{{\scriptstyle 1}}{{\scriptstyle 2}}$}a_m{%
\mathbf{Q}}_{J^2}{\mathbf{e}}_{m-1}  \notag \\
&\hspace{-.75in} = \mbox{$\frac{{\scriptstyle 1}}{{\scriptstyle 2i}}$} \big( %
\mbox{$\frac{{\scriptstyle 1}}{{\scriptstyle 2}}$}a_{m+2}{\mathbf{e}}_{m+2} +%
\mbox{$\frac{{\scriptstyle 1}}{{\scriptstyle 2}}$}a_{m+1}{\mathbf{e}}_m %
\big) - \mbox{$\frac{{\scriptstyle 1}}{{\scriptstyle 2i}}$} \big( %
\mbox{$\frac{{\scriptstyle 1}}{{\scriptstyle 2}}$}a_m{\mathbf{e}}_m +%
\mbox{$\frac{{\scriptstyle 1}}{{\scriptstyle 2}}$}a_{m-1}{\mathbf{e}}_{m-2} %
\big)  \notag \\
&\hspace{-.25in} -\mbox{$\frac{{\scriptstyle 1}}{{\scriptstyle 2}}$} \big( %
\mbox{$\frac{{\scriptstyle 1}}{{\scriptstyle 2i}}$}a_{m+2}{\mathbf{e}}_{m+2}
-\mbox{$\frac{{\scriptstyle 1}}{{\scriptstyle 2i}}$}a_{m+1}{\mathbf{e}}_m %
\big) - \mbox{$\frac{{\scriptstyle 1}}{{\scriptstyle 2}}$} \big( %
\mbox{$\frac{{\scriptstyle 1}}{{\scriptstyle 2i}}$}a_m{\mathbf{e}}_m -%
\mbox{$\frac{{\scriptstyle 1}}{{\scriptstyle 2i}}$}a_{m-1}{\mathbf{e}}_{m-2} %
\big)  \notag \\
&\hspace{-.75in} = \mbox{$\frac{{\scriptstyle 1}}{{\scriptstyle 2i}}$}%
(a^2_{m+1}-a^2_m) {\mathbf{e}}_m =%
\mbox{$\frac{{\scriptstyle
1}}{{\scriptstyle 2i}}$}(-2m{\hbar}^2) {\mathbf{e}}_m = i\hbar \, {\mathbf{Q}%
}_{J^3}{\mathbf{e}}_m.  \notag
\end{align}
Similarly 
\begin{equation*}
[{\mathbf{Q}}_{J^2}, {\mathbf{Q}}_{J^3}]{\mathbf{e}}_m = i\hbar \, %
\mbox{$\frac{{\scriptstyle 1}}{{\scriptstyle 2}}$} \big( a_{m+1}{\mathbf{e}}%
_{m+1} + a_m{\mathbf{e}}_{m-1} \big) = i\hbar \, {\mathbf{Q}}_{J^1}{\mathbf{e%
}}_m
\end{equation*}
and 
\begin{equation*}
[{\mathbf{Q}}_{J^1}, {\mathbf{Q}}_{J^3}]{\mathbf{e}}_m = i\hbar \, %
\mbox{$\frac{{\scriptstyle 1}}{{\scriptstyle 2}}$} \big( a_{m+1}{\mathbf{e}}%
_{m+1} + a_m{\mathbf{e}}_{m-1} \big) = -i\hbar \, {\mathbf{Q}}_{J^2}{\mathbf{%
e}}_m.
\end{equation*}
The operators $\mbox{$\frac{{\scriptstyle 1}}{{\scriptstyle i\hbar }}$}\,%
\boldsymbol{Q}_{J^{1}}$, $%
\mbox{$\frac{{\scriptstyle 1}}{{\scriptstyle
i\hbar}}$}{\mathbf{Q}}_{J^2} $, $%
\mbox{$\frac{{\scriptstyle
1}}{{\scriptstyle i\hbar }}$}\,\boldsymbol{Q}_{J^{3}}$ are skew symmetric
and satisfy the commutation relations of the generators of the Lie algebra $%
\mathop{\mathrm{so}}\nolimits (3)$, namely 
\begin{equation*}
\lbrack \mbox{$\frac{{\scriptstyle 1}}{{\scriptstyle i\hbar }}$}\,%
\boldsymbol{Q}_{J^{1}},%
\mbox{$\frac{{\scriptstyle 1}}{{\scriptstyle i\hbar
}}$}\,\boldsymbol{Q}_{J^{2}}]=%
\mbox{$\frac{{\scriptstyle 1}}{{\scriptstyle
i\hbar }}$}\,\boldsymbol{Q}_{J^{3}},\,\,[%
\mbox{$\frac{{\scriptstyle 1}}{{\scriptstyle i\hbar
}}$}\,\boldsymbol{Q}_{J^{2}},%
\mbox{$\frac{{\scriptstyle 1}}{{\scriptstyle
i\hbar }}$}\,\boldsymbol{Q}_{J^{3}}]=%
\mbox{$\frac{{\scriptstyle
1}}{{\scriptstyle i\hbar }}$}\,\boldsymbol{Q}_{J^{1}},\,\,[%
\mbox{$\frac{{\scriptstyle 1}}{{\scriptstyle i\hbar }}$}\,\boldsymbol{Q}%
_{J^{3}},\mbox{$\frac{{\scriptstyle 1}}{{\scriptstyle i\hbar }}$}\,%
\boldsymbol{Q}_{J^{1}}]=%
\mbox{$\frac{{\scriptstyle 1}}{{\scriptstyle i\hbar
}}$}\,\boldsymbol{Q}_{J^{2}}.
\end{equation*}%
The representation space $\mathfrak{H}$ of $\mathop{\mathrm{so}}\nolimits(3)$
has dimension $2s+1$. Thus we have constructed a (2s+1)-dimensional
representation on $\mathfrak{H}$ of the Lie algebra $\mathop{\mathrm{so}}%
\nolimits (3)$. This representation gives rise to the Lie algebra
homomorphism 
\begin{equation}
\rho : \mathop{\mathrm{so}}\nolimits (3) \rightarrow \mathop{\mathrm{gl}}%
\nolimits (\mathfrak{H}, \mathbb{R} ): \widehat{J^i\, e_i} \mapsto {\mathbf{Q%
}}_{J^i}.  \label{eq-sec6ss2onenw}
\end{equation}
Below we show that the map $\rho $ (\ref{eq-sec6ss2onenw}) can be integrated
to the Lie group homomorphism 
\begin{equation}
R: \mathop{\mathrm{SO}}\nolimits (3) \rightarrow \mathop{\mathrm{Gl}}%
\nolimits (\mathfrak{H}, \mathbb{R} ): g \mapsto {\mathop{\mathrm{Ad}}%
\nolimits}_{\exp \rho (\log g)}.  \label{eq-sec6ss2onestar}
\end{equation}
Because the representation of $\mathop{\mathrm{so}}\nolimits (3)$ on $%
\mathfrak{H}$ is irreducible, it follows that the representation of $%
\mathop{\mathrm{SO}}\nolimits (3)$, given by 
\begin{equation*}
R(g): \mathfrak{H} \rightarrow \mathfrak{H}:e_{\mathbf{m}} \mapsto R(g)e_{%
\mathbf{m}}, \quad \text{for every\, \, }g\in G
\end{equation*}
is irreducible and corresponds to spin $s \in \mathbb{N}$. \medskip

In general, let $\mathfrak{g}$ be the linear Lie algebra of the linear Lie
group $\mathfrak{G}$. Let 
\begin{equation*}
\rho : \mathfrak{g} \rightarrow \mathop{\mathrm{gl}}\nolimits (V, \mathbb{R}
): X \mapsto \rho (X)
\end{equation*}
be the Lie algebra homomorphism associated to the representation of $%
\mathfrak{g}$ on the finite dimensional real vector space $V$ given by $%
\rho(X):V \rightarrow V$, for every $X \in \mathfrak{g}$. \medskip

\noindent \textbf{Claim 6.2.1} Let 
\begin{equation}
R:\mathfrak{G}\rightarrow \mathop{\mathrm{Gl}}\nolimits(V,\mathbb{R}):\exp
X\mapsto {\mathop{\mathrm{Ad}}\nolimits}_{\exp \rho (X)}.
\label{eq-sec6ss2onedagger}
\end{equation}%
Then the map $R$ is a local Lie group homomorphism, which is defined in an
open neighborhood $\mathcal{U}$ of the identity element of $\mathfrak{G}$
where $\exp :\mathcal{V}\subseteq \mathfrak{g}\rightarrow \mathcal{U}%
\subseteq \mathfrak{G}$ is invertible. \medskip

\noindent \textbf{Proof}. For $X$, $Y$, and $Z \in \mathcal{V}$ applying the
Lie algebra homomorphism $\rho $ to the Cambell-Baker-Hausdorff formula \cite%
{dynkin} 
\begin{align}
Z(X,Y) & = \log (\exp X \, \exp Y)  \notag \\
& =\sum_{n >0} \frac{(-1)^{n-1}}{n} \sum_{\overset{r_i+s_i >0}{1\le i \le n}%
} \frac{\big( \sum^n_{i=1} (r_i+s_i)\big)^{-1}}{r_1!s_1! \cdots r_n!s_n!} \,
T_{r_1s_1 \cdots r_ns_n},  \label{eq-sec6ss2dagger}
\end{align}
where 
\begin{align}
T_{r_1s_1 \cdots r_ns_n} & =\left\{ 
\begin{array}{l}
\underbrace{\lbrack X , \lbrack X, \ldots, \lbrack X}_{r_1}, \underbrace{%
\lbrack Y, \ldots , \lbrack Y,}_{s_1} \cdots \underbrace{\lbrack X , \lbrack
X, \ldots, \lbrack X}_{r_n}, \underbrace{\lbrack Y, \ldots , \lbrack Y, Y ]}%
_{s_n} \rbrack \cdots \rbrack , \\ 
\rule{0pt}{16pt} 0, \quad \mbox{if $s_n >1$ or if $s_n=0$ and $r_n> 1$}%
\end{array}
\right.  \notag \\
\rule{0pt}{16pt} & = X+Y + \mbox{$\frac{{\scriptstyle 1}}{{\scriptstyle 2}}$}
\lbrack X,Y] + \mbox{$\frac{{\scriptstyle 1}}{{\scriptstyle 12}}$} \lbrack
X, [X,Y]] - \mbox{$\frac{{\scriptstyle 1}}{{\scriptstyle 12}}$}\lbrack Y,
[X,Y]] + \cdots  \notag
\end{align}
shows that 
\begin{equation}
\rho (Z) (\rho (X), \, \rho (Y)) = \log \, \exp \rho (X)\, \exp \rho (Y).
\label{eq-sec6ss2twonew}
\end{equation}
Therefore 
\begin{align}
R(\exp X\, \exp Y) & = R( \exp Z) = {\mathop{\mathrm{Ad}}\nolimits}_{\exp
\rho (Z)}, \text{\quad by definition}  \notag \\
& = {\mathop{\mathrm{Ad}}\nolimits}_{\exp \rho (X) \, \exp \rho (Y)}, \text{%
\quad by (\ref{eq-sec6ss2twonew})}  \notag \\
& = {Ad}_{\exp \rho (X)}\, {\mathop{\mathrm{Ad}}\nolimits}_{\exp \rho (Y)} =
R(\exp X)\, R(\exp Y),  \notag
\end{align}
that is, $R$ is a local group homomorphism. \hfill $\square $ \medskip

\noindent \textbf{Corollary 6.2.2} If $\mathfrak{G}$ is compact, then the
map $R$ (\ref{eq-sec6ss2onestar}) is a homomorphism of Lie groups. \medskip

\noindent \textbf{Proof}. Because $\mathfrak{G}$ is compact and the map $R$
is continuous, the exponential maps $\exp : \mathfrak{g} \rightarrow 
\mathfrak{G}$ and $\exp : \rho (\mathfrak{g}) \rightarrow R(\mathfrak{G})$
are surjective. Using the preceding observation, the corollary follows.
\hfill $\square $

\end{document}